\documentclass{tut2016}

\usepackage[colorlinks=true,breaklinks=true,bookmarks=false,urlcolor=blue,%
  citecolor=blue,linkcolor=blue,bookmarksopen=false,draft=false]{hyperref}

\usepackage{amsfonts, mathtools, tcolorbox}
\usepackage{dsfont, enumitem}
\usepackage{subfigure, epstopdf} 
\usepackage{algorithm, algorithmic}

\newcommand{\R}{{\mathbb{R}}}
\renewcommand{\Pr}{{P}}
\newcommand{\E}{{E}}
\newcommand{\xsp}{{\mathcal{X}}}

\newcommand{\risk}{{\mathcal{R}}}
\newcommand{\drorisk}{{\mathcal{R}_{\delta}}}
\newcommand{\ermin}{{\hat{\theta}_n^{\,\textnormal{erm}}}}
\newcommand{\drmin}{{\hat{\theta}_n^{\,\textnormal{dro}}}}
\newcommand{\Let}{\coloneqq}

\newcommand{\T}{{\intercal}}

\newcommand{\dconv}{{\overset{dist.}{\longrightarrow}}}

\begin{document}
\CHAPTERNO{}
\TITLE{Statistical Analysis of Wasserstein Distributionally Robust Estimators}    

\AUBLOCK{%
  \AUTHOR{Jose Blanchet}
  \AFF{Management Science \& Engineering, Stanford University, CA 94305, 
       \EMAIL{jose.blanchet@stanford.edu}} 
  \AUTHOR{Karthyek Murthy}
  \AFF{Singapore University of Technology and Design, Singapore 487372, 
       \EMAIL{karthyek\_murthy@sutd.edu.sg}} 
    \AUTHOR{Viet Anh Nguyen}
  \AFF{Management Science \& Engineering, Stanford University, CA 94305 and VinAI Research, Vietnam, 
       \EMAIL{viet-anh.nguyen@stanford.edu}}
} 

\CHAPTERHEAD{Blanchet, Murthy and Nguyen: Statistical Analysis of Wasserstein Distributionally Robust Estimators}

\ABSTRACT{%
We consider statistical methods which invoke a min-max distributionally
robust formulation to extract good out-of-sample performance in data-driven
optimization and learning problems. Acknowledging the distributional
uncertainty in learning from limited samples, the min-max formulations
introduce an adversarial inner player to explore unseen covariate data. The
resulting Distributionally Robust Optimization (DRO) formulations, which
include Wasserstein DRO formulations (our main focus), are specified using optimal transportation phenomena. Upon describing how these
infinite-dimensional min-max problems can be approached via a
finite-dimensional dual reformulation, the tutorial moves into its main
component, namely, explaining a generic recipe for optimally selecting the
size of the adversary's budget. This is achieved by studying the limit behavior of an optimal transport projection formulation arising from an inquiry on the
smallest confidence region that includes the unknown population risk
minimizer. Incidentally, this systematic prescription coincides with those
in specific examples in high-dimensional statistics and results in error
bounds that are free from the curse of dimensions. Equipped with this
prescription, we present a central limit theorem for the DRO estimator and
provide a recipe for constructing compatible confidence regions that are useful for uncertainty quantification. The rest of
the tutorial is devoted to insights into the nature of the optimizers selected by
the min-max formulations and additional applications of optimal transport
projections.
} 


\maketitle  

\section{Introduction}
\label{sec:intro}

Data-driven decision making has permeated virtually every aspect of
operations research (OR) and management science (MS). This proliferation has been made possible
thanks to our increased ability to collect a gigantic amount of data and the
availability of computational resources which enable the solution of complex
uncertainty-informed optimization problems. In many OR/MS tasks, decisions
are chosen as prescriptions to improve future performance. Hence, it is
imperative to recognize that the available evidence (often based on the limited data collected from previous experience or presumably similar environments) might deflect from the future environment in which the decision will be applied. This recognition ignites the field of decision making under uncertainty, and an emerging framework for robust decision making and statistical analysis is that of distributionally robust optimization (DRO).

While DRO formulations have received substantial
attention in the OR/MS community during the last decades \cite{ref:bental2013robust, ref:wiesemann2014distributionally, ref:shapiro2017distributionally, ref:chen2020robust, ref:vanparys2021from, ref:pflug2007ambiguity, ref:pflug20121/N}, recent years have witnessed a significant amount of interest in statistical properties of data-driven DRO-based decision rules and
associated inference obtained by these of decisions. The main goal of this tutorial is to discuss some statistical properties enjoyed by decisions obtained from DRO formulations as well as the associated techniques that are used to analyze these types of
decisions. Our focus is on DRO formulations for which the uncertainty
region is described in terms of so-called optimal transport discrepancies (which include the Wasserstein\footnote{The correct spelling of Wasserstein appears to be Vasershtein (in honor of L.~N.~Vasershtein), but most of the literature uses the Wasserstein spelling, so we keep this spelling in this tutorial.} distances as special cases). 

The literature that connects optimal transport and data analytics has grown substantially in recent years. The text~\cite{ref:peyre2019computational} provides many examples in which optimal transport is used in areas such as computer vision and machine learning, with an emphasis on computational methods. The Wasserstein distance plays a key role in the design of a popular generative artificial intelligence algorithm known as the Wasserstein Generative Adversarial Network~\cite{ref:arjovsky2017wasserstein}. Another popular application of optimal transport in machine learning relates to study of adversarial attacks \cite{ref:goodfellow2015explaining, sinha2018certifiable, ref:madry2018towards}. Other applications include optimal transport and domain adaptation (i.e., transferring a model trained on one environment to another) \cite{ref:courty2017optimal}; optimal transport and missing data~\cite{ref:muzellec2020missing}, optimal transport in the context of Bayesian computation~\cite{ref:srivastava2015wasp, ref:bernton2019approximate}; and optimal transport in deconvolution and denoising \cite{ref:rigollet2018entropic, ref:nguyen2019bridging}, among others. In addition, there are burgeoning variants of the optimal transport distance, including the unbalanced optimal transport~\cite{ref:chizat2018unbalanced}, subspace robust Wasserstein distance~\cite{ref:paty2019subspace}, sliced Wasserstein distance~\cite{ref:kolouri2019generalized}, tree-sliced Wasserstein distance~\cite{ref:le2019tree}, etc. While these applications and variants of the optimal transport distance are of great interest, there is also a rich statistical structure underlying the use of optimal transport in these settings. The focus of this tutorial is on data-driven DRO formulations which utilize the optimal transport theory to inform distributional shifts from the empirical distribution. We believe, however, that the techniques that we discuss, including the projection analysis in the Wasserstein geometry and the hypothesis testing tools, could be extended to provide statistical insights in many of the applications discussed above.  

To facilitate our discussion of the statistical properties of Wasserstein-DRO estimators, we first discuss in Sections~\ref{intro-elements} and \ref{intro-dro} some basic definitions.

\subsection{Introductory elements\label{intro-elements}}
In order to quickly go to the heart of our technical discussion, let us
introduce a generic expected loss minimization, which can be seen as an
idealized decision making problem under full information. In particular,
given the loss $\ell :\mathcal{X}\times \Theta \rightarrow 
\mathbb{R}$, suppose that we wish to solve 
\begin{align}
\inf_{\theta \in \Theta }E_{\Pr_{\ast }}\left[ \ell (X,\theta )%
\right] ,
\label{prm}
\end{align} %
where $X$ is a random vector taking values in $\mathcal{X}\subseteq \mathbb{R%
}^{m}$, $\Theta \subseteq \mathbb{R}^{d}$, and $%
\Pr_{\ast }$ denotes the unknown distribution of $X$. The decision space in
this case is $\Theta $. For ease of notation, let 
\begin{equation*}
\mathcal{R}(\Pr ,\theta )\coloneqq E_{\Pr }\left[ \ell
(X,\theta )\right] 
\end{equation*}%
denote the expected loss (or risk) associated with the parameter/decision
choice $\theta $ when evaluated under the distributional assumption $X\sim
\Pr $. Then (\ref{prm}) gets equivalently written as 
\begin{equation*}
\mathcal{R}_{\ast }(P_{\ast })\coloneqq\inf_{\theta \in \Theta }%
\mathcal{R}(\Pr_{\ast },\theta ),
\end{equation*}%
and $\mathcal{R}_{\ast }(P_{\ast })$ denotes the optimal risk.

Let $(X_1,\ldots,X_n)$ be an independent and identically distributed (i.i.d.) sample from the unknown distribution $\Pr_\ast$ (to ease notation we will write $\E[f(X_1,...,X_n)]$ for any $f$ with well defined expectation, that is, we will not write subscripts with the $n$-th fold product of $\Pr_\ast$). A standard approach towards solving \eqref{prm} entails minimizing the \textit{%
empirical risk}, 
\begin{align}
\inf_{\theta \in \Theta }~\mathcal{R} (\Pr_n,\theta),
\label{erm}
\end{align}
where the empirical distribution $\Pr_n \coloneqq   \frac{1}{n}%
\sum_{i=1}^n\delta_{X_i}$ is plugged-in place of the unknown distribution $%
\Pr_\ast$ in the population risk $\mathcal{R}(\Pr_\ast,\theta)$ in %
\eqref{prm}. This is indeed natural as the sample average loss, $\mathcal{R} %
(P_n,\theta)$, constitutes an unbiased (that is, $\E[\risk(\Pr_n,\theta)] = \risk(\Pr_\ast, \theta)$), minimum-variance estimator of $\mathcal{R} (\Pr_\ast,\theta)$, for
any $\theta \in \Theta $ (see, for example, \cite[Chapter 5]{shapiro2014lectures}).

For a solution to \eqref{erm}, namely $\hat{\theta}_n^{\,%
\textnormal{erm}}  \in \text{arg}\,\text{min} _{\theta \in
\Theta} \mathcal{R} (P_n,\theta)$, we can however only assert the
out-of-sample risk to be witnessed with the empirical optimum $\hat{\theta}%
_n^{\,\textnormal{erm}} $ always exceeds the expected in-sample risk:
indeed, 
\begin{align*}
\mathcal{R} \big(\Pr_\ast, {\hat{\theta}_n^{\,\textnormal{erm}} %
} \big) \geq \inf_{\theta \in \Theta} \mathcal{R} (\Pr_\ast,
\theta) = \inf_{\theta \in \Theta} E \left[ \mathcal{R} %
(\Pr_n,\theta)\right] \geq E \left[\inf_{\theta \in \Theta} \mathcal{R%
} \big(\Pr_n, \theta\big)\right] = E\left[\mathcal{R}  \big(%
P_n,\hat{\theta}_n^{\,\textnormal{erm}}  \big) \right],
\end{align*}
highlighting that the in-sample optimal risk suffers from an optimistic
bias. The gap $\mathcal{R} \big(\Pr_\ast, \hat{\theta}_n^{\,%
\textnormal{erm}} \big) - \mathcal{R}  \big(P_n,\hat{\theta}%
_n^{\,\textnormal{erm}}  \big)$, which quantifies the post-decision
disappointment, can often be large, and remarkably so in high-dimensional
settings. This phenomenon gets referred to as ``optimizer's curse" or
``overfitting", based on the context; see \cite{ref:kuhn2019wasserstein} and
references therein for a more detailed discussion.
We next explore the DRO approach for mitigating this difficulty. 

\subsection{Distributionally robust optimization formulations.}\label{intro-dro} 

A recent approach which has gained prominence in mitigating optimistic bias and other considerations discussed earlier is a \textit{%
distributionally robust} variant of \eqref{erm} which accounts for the
distributional uncertainty in utilizing the empirical measure $\Pr_n$ as a
proxy for $\Pr_\ast$. The effect of this distributional uncertainty is
incorporated by instead minimizing the worst-case risk, 
\begin{align}
\mathcal{R}_{\delta} (\Pr_n,\theta) \ \coloneqq \sup_{\Pr \, \in \,%
\mathcal{U}_\delta(\Pr_n)} \mathcal{R} (\Pr,\theta),  \label{dr-risk}
\end{align}
evaluated over a set $\mathcal{U}_\delta(\Pr_n)$ of probability
distributions which are \textit{plausible as a candidate for $\Pr^\ast$ in
the task of solving \eqref{prm}}. The set $\mathcal{U}_\delta(\Pr_n)$ is
referred as the \textit{distributional ambiguity set}. The resulting
optimization problem, 
\begin{align}
\inf_{\theta \in \Theta } \mathcal{R}_{\delta} (P_n , \theta)
= \inf_{\theta \in \Theta } \sup_{\Pr \, \in \,\mathcal{U}%
_\delta(\Pr_n)} E _{\Pr} \left[ \ell(X,\theta)\right],  \label{dro}
\end{align}
is referred as a DRO formulation for
solving \eqref{prm}. One may view the inner supremum as the effect of an
adversary free to explore the implications of varying the benchmark model,
which is $P_n$ in this case, within the ambiguity set $\mathcal{U}%
_\delta(P_n)$. The DRO formulation then seeks a choice, denoted by $\hat{%
\theta}_n^{\,\textnormal{dro}}  \in \Theta$, which minimizes the
worst-case expected cost in \eqref{dro} induced by the adversary specified
by the distributional ambiguity $\mathcal{U}_\delta(P_n)$. 

We now introduce the notion of \textit{optimal transport costs} between probability distributions. Let $\mathcal{P%
}(\mathcal{X} )$ denote the collection of probability measures
defined on the Borel space of $\mathcal{X}$. The space $\mathcal{X}$ is assumed to be a complete separable metric space.  For the purpose of this tutorial one might think of $\mathcal{X}$ as the Euclidean space. The reader is referred to \cite{ref:villani2003topics, ref:villani2008optimal, ref:santambrogio15optimal} for an introduction to optimal transport theory.

\begin{definition}[Optimal transport costs, Wasserstein distances]
\textnormal{Given a lower semicontinuous function $c:\mathcal{X}  %
\times \mathcal{X}  \rightarrow [0,\infty]$, the optimal transport
cost $D_c(P,Q)$ between any two distributions $P,Q \in \mathcal{P}(\mathcal{X%
} )$ is defined as 
\begin{align*}
D_c(P,Q) = \min_{\pi \in \Pi(P,Q)} E_{\pi}\left[c(X,X^\prime)\right]
\end{align*}
where $\Pi(P,Q)$ denotes the set of all joint distributions of the random
vector $(X,X^\prime)$ with marginal distributions $P$ and $Q$, respectively.
If we specifically take $c(x,x^\prime) = \Vert x - x^\prime \Vert^r$, for $r
\in [1,\infty)$, we obtain a Wasserstein distance of type $r$ by letting $%
W_r(P,Q) = \left\{D_c(P,Q)\right\}^{1/r}$. }  \label{defn:WD}
\end{definition}

The quantity $D_c(P_{n},P)$ may be interpreted as the cheapest way to
transport mass from the distribution $P_{n}$ to the mass of another
probability distribution $P$, while measuring the cost of transportation
from location $x \in \mathcal{X} $ to location $y \in \mathcal{X} %
$ in terms of the transportation cost $c(x,y)$. Equipped with a
notion of distance between probability distributions, a natural formulation
of the distributional ambiguity is given by 
\begin{align}
\mathcal{U}_\delta(\Pr_n) = \left\{ \Pr: D_c (\Pr_n, \Pr) \leq \delta
\right\},  \label{u-delta-defn}
\end{align}
for a suitable radius (or) budget of ambiguity, captured by the parameter $%
\delta > 0$. The Wasserstein distances $W_r(\cdot)$ serve as the canonical
choice for informing the distance $D_c$ in \eqref{u-delta-defn}. The
resulting ambiguity set, as we shall see, includes all feasible random perturbations of the form $%
\{X_i + \Delta_i: i = 1,\ldots,n\}$ to the training samples such that the
perturbations are constrained in the $L^r$ norm.  The goal of the Wasserstein
distanced based DRO procedure is choosing a decision that also hedges
against these adversarial perturbations, thus introducing adversarial
robustness into settings where the quality of optimal solutions is
sensitive to incorrect model assumptions. As argued in \cite{ref:kuhn2019wasserstein}, the DRO formulation of the type \eqref{dro} can be motivated from an axiomatic approach; see \cite{ref:delage2019dice} and \cite{GILBOA1989141}. The Wasserstein-DRO formulation clearly explores the impact of out-of-sample scenarios as explained in~\cite{ref:esfahani2018data}; other forms of uncertainty sets include divergence measures~\cite{ref:bental2013robust, ref:bayraksan2015data} and moment-based constraints~\cite{ref:delage2010distributionally, ref:goh2010distributionally}. The work \cite{ ref:rahimian2019distributionally} provides a comprehensive review of DRO methods with a special emphasis on optimization techniques and results. Our focus here is on statistical output analysis tasks encompassing optimal selection of the  uncertainty size $\delta$ and characterizing associated asymptotic normality and confidence regions.   

We also note that a modeler may choose to depart from the choice $%
c(x,x^\prime) = \Vert x - x^\prime \Vert^r$ in problems with differing
geometries. Though we shall be restricting attention to this canonical choice in most
examples here for the sake of simplicity, one may view the transportation
cost $c(\cdot)$ as a powerful modeling tool in exploring the impact of
distributional uncertainty. The examples in \cite%
{ref:blanchet2020distributionally,ref:blanchet2017data} serve to illustrate
the improved out-of-sample performances one may obtain by suitably
incorporating the geometry of the problem in informing the optimal transport
costs.

\subsection{Organization} The rest of the tutorial is organized as follows. In Section \ref{sec:refor} we briefly review some duality results and examples which are used to motivate statistical properties of Wasserstein-DRO estimators. We then move on to discuss how to optimally select the size of the ambiguity set in Wasserstein-DRO estimators using a certain hypothesis testing criterion which is connected to projection techniques. This is done in Section \ref{sec:WPF-Intro-Dev}. We provide a discussion of alternative methods, including cross validation and finite sample guarantees. The optimal ambiguity set size is further studied from a confidence region perspective together with asymptotic normality under the optimal choice and finite sample guarantees are discussed in Section \ref{sec:dro-clt}. Final considerations and conclusions are given in our last section, namely, Section~\ref{sec:conclusions}.

\section{Dual reformulation and examples}
\label{sec:refor} 

Solving the DRO formulation~\eqref{dro} naturally requires
evaluating the worst-case risk $\mathcal{R}_{\delta} (P_n,\theta)$
for any $\theta \in \Theta $. From Definition \ref{defn:WD} and from the
formulation of $\mathcal{U}_\delta(P_n)$ in \eqref{dr-risk}, we have 
\begin{align}
\mathcal{R}_{\delta} (P_n,\theta) = \sup \left\{ E_\pi\left [
\ell(X^\prime, \theta)\right]: \pi \in \Pi(\Pr_n,\Pr), \ E _\pi\left[
c(X,X^\prime)\right] \leq \delta \right\}.  \label{primal-reform}
\end{align}
Observe that the joint measure $\pi \in \Pi(P_n, P)$ can be written in terms of marginal
constraints of the form $\pi (A \times \mathcal{X} ) = \Pr_n(A)$ or $%
\pi (\mathcal{X}  \times A) = \Pr(A)$, for Borel subsets $A$. Though
infinitely many, these constraints are linear over the measure $\pi$. Thus
the objective and the constraints in the evaluation of $\mathcal{R}_{\delta} %
(P_n,\theta)$ are linear in the variable $\pi$. Thanks to this linear
programming structure, one can reformulate this infinite-dimensional maximization problem using duality theory~\cite{ref:esfahani2018data, ref:blanchet2019quantifying, ref:gao2016distributionally, ref:zhao2018data}. 

\begin{theorem}[Strong duality]
\label{thm:duality} Suppose the transportation cost $c:\mathcal{X}  %
\times \mathcal{X}  \rightarrow [0,\infty]$ satisfies $c(x,x) = 0$
for all $x \in \mathcal{X} $. Then for any reference probability
distribution $\Pr_{\mathnormal{ref}}$ and upper semicontinuous $f: \mathcal{X%
} \rightarrow \mathbb{R} $ satisfying $E _{\Pr_{\mathnormal{ref%
}}}\vert f(X) \vert < \infty$, we have 
\begin{align}
\sup_{\Pr: D_c(\Pr_{\mathnormal{ref}},\Pr) \leq \delta} E _{\Pr} %
\left[ f(X)\right] = \inf_{\lambda \geq 0} \ \ \lambda \delta + E %
_{\Pr_{\mathnormal{ref}}}\left[f_{\lambda}(X)\right] ,  \label{dual-reform}
\end{align}
where $f_{\lambda}(x) \Let \sup_{z \in \mathcal{X} } \left\{
f(z) - \lambda c(x,z)\right\}$.
\end{theorem}
\noindent 
Notice that Theorem~\ref{thm:duality} holds for any reference measure $P_{ref}$, which also includes the case of interest in this tutorial of the empirical measure $P_n$. Moreover, the restriction of $c(x,x)=0$ and integrability can be further relaxed; see~\cite{ref:kent2021frank-wolfe}. 

It is instructive to verify the duality in the case of $\mathcal{X} $
being a finite set. Suppose that $\mathcal{X}  = \{x_1,\ldots,x_k\}$
and $\Pr_{\mathnormal{ref}} (X = x_i) = p_i$ for $i = 1,\ldots,k$; then 
\begin{align*}
\sup_{\Pr: D_c(\Pr_{\mathnormal{ref}},\Pr) \leq \delta} E _{\Pr} %
\left[ f(X)\right] = \left\{ 
\begin{array}{cl}
\max & \sum_{j}f(x_j)\pi_{i,j} \\ 
\mathrm{s.t.} & \pi \in \mathbb{R} _+^{k \times k} \\ 
& \sum_{j=1}^k \pi_{i,j} = p_i \quad \forall i =1, \ldots, k \\ 
& \sum_{i,j=1}^k c(x_i,x_j)\pi_{i,j} \leq \delta%
\end{array}
\right.
\end{align*}
is a finite-dimensional linear program. In this case, the dual linear program is expressly written as 
\begin{align*}
\min \left\{ \lambda \delta + \sum_{i=1}^k p_i \nu_i: \lambda \geq 0, \nu_i
\geq f(x_j) - \lambda c(x_i,x_j)~~\forall(i,j)\in \{1,\ldots,k\}^2 \right\},
\end{align*}
which equals the RHS in \eqref{dual-reform}.

\begin{corollary}
Suppose that $\ell(x,\theta)$ is upper semicontinuous in $x$, for any $\theta \in \Theta$. Then for the transportation cost $c(x,x^\prime) = \Vert
x - x^\prime \Vert^r$, where $r \in [1,\infty)$, we have for any $\theta \in
\Theta$ 
\begin{align*}
\mathcal{R}_{\delta} (P_n,\theta) = \inf_{\lambda \geq 0} \ \ \lambda
\delta + \frac{1}{n} \sum_{i=1}^n \sup_{\Delta : X_i + \Delta \in \mathcal{X}
} \big\{ \ell(x_i + \Delta, \theta) - \lambda \Vert \Delta\Vert^r %
\big\}.
\end{align*}
\label{cor:duality}
\end{corollary}

\begin{remark}[Structure of the adversarial distribution attaining the
sup in \eqref{dual-reform}]
\textnormal{An optimal coupling $\pi^\ast$ that attains the maximum in~\eqref{primal-reform}, if it exists, can be written as the joint law of $%
(X,X_{\mathnormal{adv}})$ satisfying
\begin{align*}
&X_{\mathnormal{adv}}  \in \arg\max_{z \in\mathcal{X}} \left\{
\ell(z,\theta) - \lambda c(x,X)\right\} \text{ a.s.}, \\
&E _{\pi^\ast}[c(X,X_{\mathnormal{adv}} )] = \delta, \text{
and } \\
&E _{\pi^\ast}[f(X_{\mathnormal{adv}} )] = \sup_{\Pr: D_c(\Pr_{%
\mathnormal{ref}},\Pr) \leq \delta} E _{\Pr} \left[ f(X)\right],
\end{align*}
where the minimization in its respective dual reformulation is attained at
some $\lambda > 0$ (see \cite[Theorem 1]{ref:blanchet2019quantifying}). Further discussion regarding the existence and the structure of the optimal coupling $\pi^\ast$ can be found in~\cite{ref:yue2021on}. In addition, for certain loss functions $\ell$ and for $P_\mathnormal{ref} = P_n$ being the empirical distribution, the locations of the atoms of $\pi^\ast$ are also known explicitly. This fact has been exploited to reformulate the distributionally robust chance-constraints~\cite{ref:chen2018data, ref:xie2021distributionally, ref:ji2021data} and to estimate the nonparametric likelihood~\cite{ref:nguyen2019optimistic}.}
\end{remark}

The strong duality result in Theorem~\ref{thm:duality} leads to tractable reformulations for various Wasserstein DRO problems. Next, we explore how this result can be applied specifically in the context of robust mean-variance portfolio allocation~\cite{ref:blanchet2020distributionally}.

\begin{example}[Worst-case expected portfolio return]
\label{eg:wc-return} Given $n$ independent samples of asset
returns $\{X_i\}_{i=1}^n$, suppose that we wish to characterize the
portfolio weights with worst-case return exceeding a target return $t;$ in
other words, we aim to identify the set 
\begin{align*}
\Theta_{\delta,t}  = \left\{\theta \in \mathbb{R} ^d:
\theta^\intercal  1 = 1, \min_{\Pr \in\, \mathcal{U}_\delta(P_n)} E %
_{\Pr}\left[ \theta^\intercal  X\right] \geq t \right\}.
\end{align*}
Taking $\ell(x, \theta) = -\theta^\intercal  x$ and the
transportation cost $c(x,x^\prime)=\Vert x - x^\prime \Vert_q^2$, for some $q \in [1,\infty)$, we evaluate the resulting worst-case risk to be 
\begin{align*}
-\min_{\Pr \in \mathcal{U}_\delta(P_n)} E _{\Pr}\left[
\theta^\intercal  X\right] &= \mathcal{R}_{\delta} (P_n,
\theta) = \inf_{\lambda \geq 0}~\lambda \delta + \frac{1}{n} \sum_{i=1}^n
\sup_{\Delta} \left\{ - \theta^\intercal  (X_i + \Delta) - \lambda
\Vert \Delta \Vert_q^2\right\}
\end{align*}
as a consequence of Corollary \ref{cor:duality}. Due to H\"{o}lder's
inequality, 
\begin{align*}
- \theta^\intercal  \Delta - \lambda \Vert \Delta \Vert_q^2 \leq
\Vert \theta \Vert_p \Vert \Delta \Vert_q - \lambda \Vert \Delta \Vert_q^2,
\qquad \Delta \in \mathbb{R} ^d,
\end{align*}
with equality attained at the choice $\Delta = \Vert \theta \Vert_p^{1 -
p/q} \mathnormal{sgn} (\theta) \vert \theta \vert^{p/q}$, with $p \in
[1,\infty)$, is such that $1/p + 1/q = 1$. Here, the sign function $%
\mathnormal{sgn} (\theta)$, the absolute value $\vert \theta \vert$,
and the exponentiation $\vert \theta \vert^{p/q}$ are applied
component-wise. Therefore, 
\begin{align*}
\mathcal{R}_{\delta} (P_n, \theta) &= -\frac{1}{n} \sum_{i=1}^n
\theta^\intercal  X_i + \inf_{\lambda > 0} \lambda \delta + \frac{1}{n%
}\sum_{i=1}^n \sup_{ \Vert \Delta \Vert_q} \left\{ \Vert \theta \Vert_p
\Vert \Delta \Vert_q - \lambda \Vert \Delta \Vert_q^2\right\} \\
&= - \theta^\intercal  E _{\Pr_n} [X] + \inf_{\lambda \geq 0}
\lambda \delta + \frac{\Vert \theta \Vert_p^2}{4\lambda} \\
&= - \theta^\intercal  E _{\Pr_n} [X] + \delta^{1/2} \Vert
\theta \Vert_p.
\end{align*}
Thus the set of portfolio weights meeting a target return is characterized
by the convex set 
\begin{align*}
\Theta_{\delta,t}  = \{\theta \in \mathbb{R}^d: \theta^\intercal %
 1 = 1,~\theta^\intercal  E _{\Pr_n}[X] \geq t +
\delta^{1/2} \Vert \theta \Vert_p \}.
\end{align*}
\end{example}

\begin{example}[Worst-case variance of return] \label{eg:variance}
Let $\mathrm{Cov}_P(X)$ be the covariance matrix of the random
vector $X$ under $P$. In a similar setting as in Example~\ref{eg:wc-return},
let $\mathcal{R} (\Pr, \theta) = \theta^\intercal  \mathrm{Cov}%
_\Pr (X) \theta$ be the variance of returns associated with a portfolio
allocation $\theta$. If the ground transportation cost is likewise taken to
be $c(x,x^\prime)=\Vert x - x^\prime \Vert_q^2$, then the worst-case
variance becomes 
\begin{equation*}
\mathcal{R}_{\delta} (P_n, \theta) = \left[\mathcal{R} (P_n,
\theta)^{1/2} + \delta^{1/2} \Vert \theta \Vert_p\right]^2. 
\end{equation*}
\end{example}

A portfolio manager can combine the worst-case mean and variance
characterizations in Examples~\ref{eg:wc-return} and~\ref{eg:variance} to
form a distributionally robust mean-variance portfolio~\cite[Theorem~1]%
{ref:blanchet2020distributionally} as below.

\begin{example}[Robust mean-variance portfolio allocation] \label{eg:drmvopt}
Suppose that we wish to construct a portfolio of risky assets
with minimum worst-case variance, while at the same time requiring a
baseline return to be met regardless of the probability distribution in $%
\mathcal{U}_\delta(\Pr_n)$. This results in the following distributionally
robust mean-variance problem
\begin{equation*}
\begin{array}{cl}
\min & \sup\limits_{P \in \mathcal{U}_\delta(P_n)}~\theta^\intercal  
\mathrm{Cov}_\Pr (X) \theta \\ 
\mathrm{s.t.} & \theta \in \mathbb{R}^d,~\theta^\intercal  1 =
1,~\min \limits_{\Pr \in\, \mathcal{U}_\delta(P_n)} E _{\Pr}\left[
\theta^\intercal  X\right] \geq t,%
\end{array}
\end{equation*}
for a specified minimum return level $t$. Then, the aforementioned
distributionally robust mean-variance problem is equivalent to a conic
optimization problem, 
\begin{equation*}
\begin{array}{cl}
\min & \left[\mathcal{R} (P_n, \theta)^{1/2} + \delta^{1/2} \Vert
\theta \Vert_p\right]^2 \\ 
\mathrm{s.t.} & \theta \in \mathbb{R}^d,~\theta^\intercal  1 =
1,~\theta^\intercal  E _{\Pr_n}[X] \geq t + \delta^{1/2} \Vert
\theta \Vert_p.%
\end{array}
\end{equation*}
Further details about the efficacy of this distributionally robust mean-variance  allocation formulation can be found in \cite{ref:blanchet2020distributionally}.
\end{example}

\begin{example}[Distributionally robust linear regression]
Given predictor-response \sloppy{\ pairs $\left\{(X_i,Y_i)
\right\}_{i = 1}^n \subset \mathbb{R} ^d \times \mathbb{R} $,}
consider the example of performing linear regression with the square loss $%
\ell(x,y,\theta) = (y - \theta^\intercal  x)^2$. Taking the
transportation cost as
\begin{align}
c\big( (x,y), (x^\prime,y^\prime)\big) = \Vert x - x^\prime \Vert_q^2 + a
\vert y - y^\prime \vert^2,  \label{trans-cost-linreg}
\end{align}
for some constant $a \in (0,\infty]$, one may similarly obtain the
corresponding worst-case loss as 
\begin{align*}
\mathcal{R}_{\delta} (P_n,\theta) &= \inf_{\lambda \geq 0} \lambda
\delta + \frac{1}{n} \sum_{i=1}^n \sup_{\Delta_x,\Delta_y} \big\{ (Y_i -
\theta^\intercal  X_i + \Delta_y - \theta^\intercal  \Delta_x
)^2 - \lambda \Vert \Delta_x \Vert_q^2 - \lambda a \Delta_y^2 \big\} \\
&= E _{\Pr_n} \left[ (Y - \theta^\intercal  X)^2\right] +
\inf_{\lambda > 0} \ \lambda \delta + \frac{E_{P_n}\left[(Y-\theta^\intercal %
 X)^2\right]}{\big(\lambda \{\Vert \theta \Vert_p^p +
a^{-p/2}\}^{-2/p}- 1 \big)^+}
\end{align*}
by using H\"{o}lder's inequality in a similar manner as in Example \ref%
{eg:wc-return}. Solving the minimization over $\lambda$ results in

\begin{align*}
\mathcal{R}_{\delta} (P_n,\theta) = \left[ E_{P_n}\big(%
Y-\theta^\intercal  X \big)^2 + \delta^{1/2} \big\{ \Vert \theta
\Vert_p^p + a^{-p/2}\big \}^{1/p}\right]^2.
\end{align*}
For the instance where $a = +\infty$, 
\begin{align}
\arg\,\min\limits_{\theta \in \Theta} \ \mathcal{R}_{\delta} %
(\Pr_n,\theta) = \arg\,\min\limits_{\theta \in \Theta} \left\{ E_{P_n}\big(%
Y-\theta^\intercal  X \big)^2 + \delta^{1/2} \Vert \theta \Vert_p
\right\},  \label{sqroot-lasso}
\end{align}
thus we recover the well-known Square-Root Lasso estimator~\cite{ref:belloni2010square} as a consequence when $q = \infty$ (and subsequently 
$p = 1$).

To contrast the qualitative behavior of the DRO estimator $\drmin \in \argmin_{\theta}~ \drorisk(\Pr_n,\theta)$  from that of the ordinary least squares solution $\ermin$, we consider the following linear regression example from \cite{ref:blanchet2021confidence}. Figures \ref{fig:scatter-p5} - \ref{fig:scatter-1} plot the realizations of $\ermin,\drmin \in \R^2$ from training independently on 1000 datasets, each of size $n = 100$, obtained by sampling from the  linear regression model $Y = \theta_\ast^\T X + \varepsilon$. The specific parameter choices are as follows:  $\varepsilon \sim \mathcal{N}(0,1)$, $X = (X_1,X_2)$ is normally distributed with $\E[X_i] = 0$, $\text{Var}[X_i] = 1$, and $\text{Cov}(X_1,X_2) = \rho \in \{-0.95, 0 , 0.95\}$. Figure \ref{fig:scatter-p5} is obtained from the linear regression model with $\theta_\ast = (0.5,0.5)^\T$, and Figure \ref{fig:scatter-1} contains estimator realizations for samples from the instance $\theta_\ast = (1,0)^\T$. We notice that $\drmin$ realizations exhibit significantly lower variability when compared to the $\ermin$ realizations in the near-collinear instances where $\vert\rho\vert = 0.95$. This stability comes, however, at the expense of a bias (shrinkage towards the origin) exhibited in the case of $\drmin$ realizations. Similar qualitative behavior is exhibited by DRO estimators in settings beyond linear regression as well, as we shall explain in terms of variation regularization in Theorem~\ref{thm:var-reg} and the asymptotic bias exhibited in the central limit theorem in Theorem~\ref{thm:joint-lim}. 
 \begin{figure}[ptbh]
\centering
\subfigure[$\rho=0.95$]{
 \includegraphics[width=1.5in]{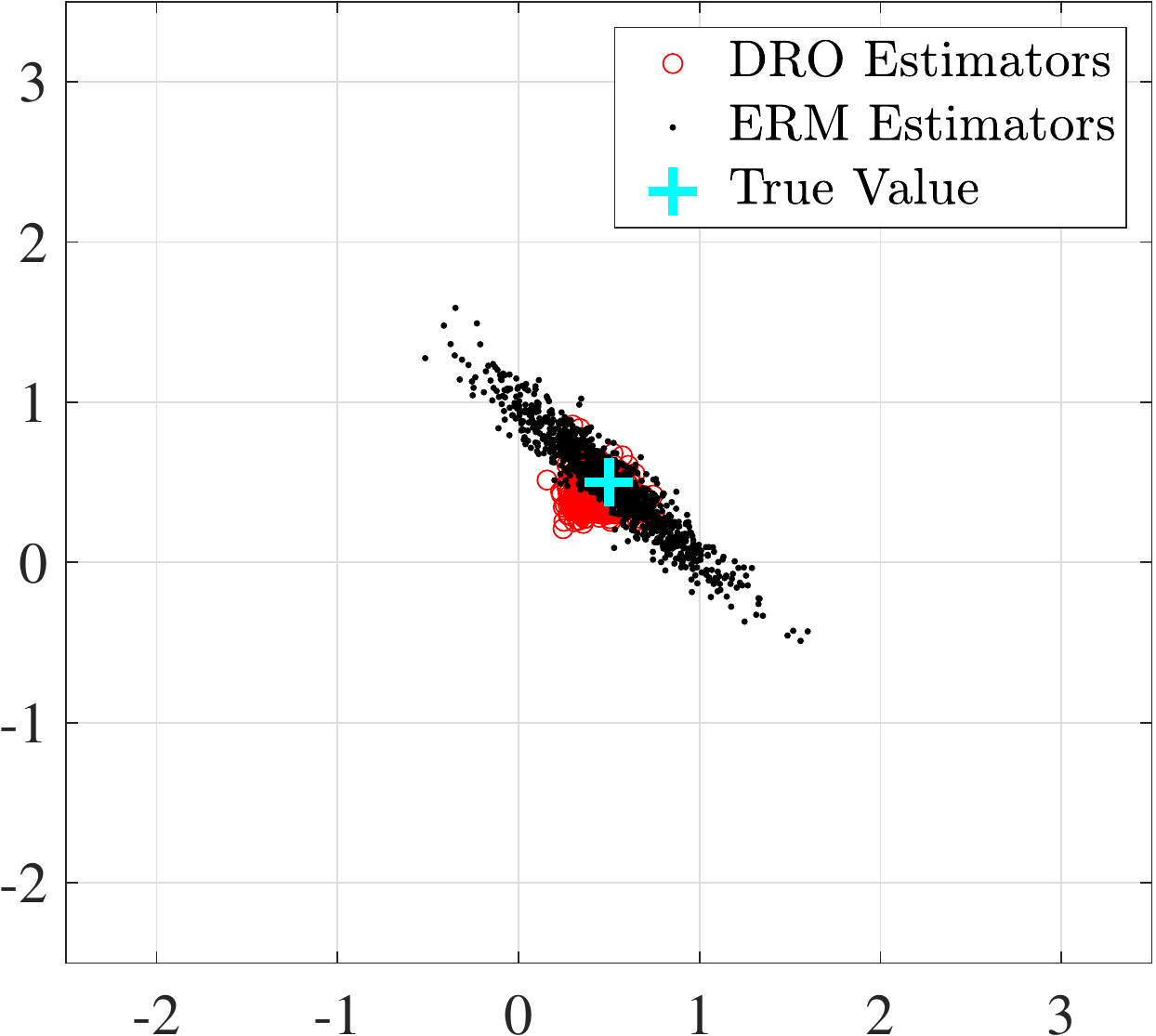}}
\subfigure[$\rho=0$]{
 \includegraphics[width=1.5in]{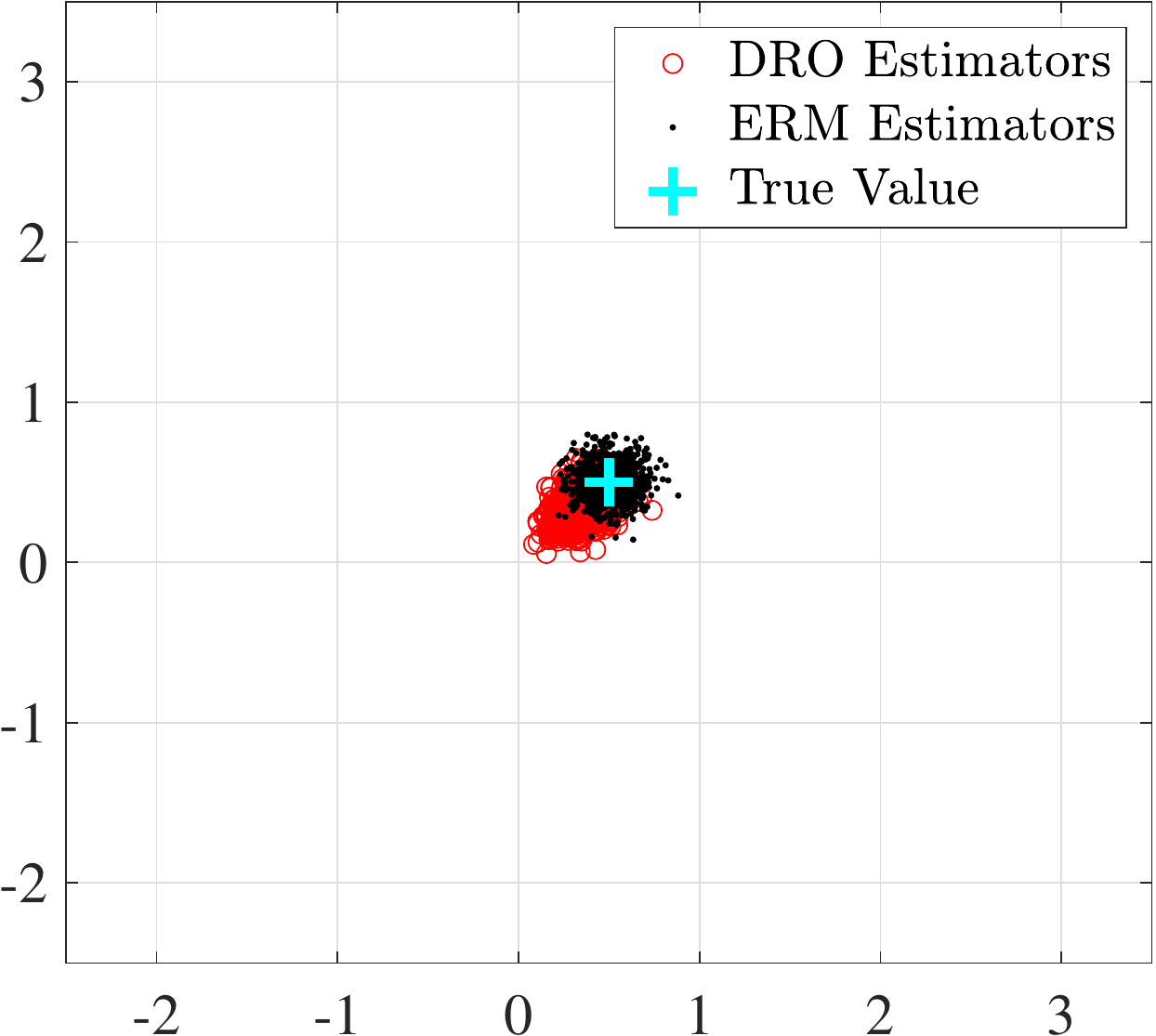}}
\subfigure[$\rho=-0.95$]{
\includegraphics[width=1.5in]{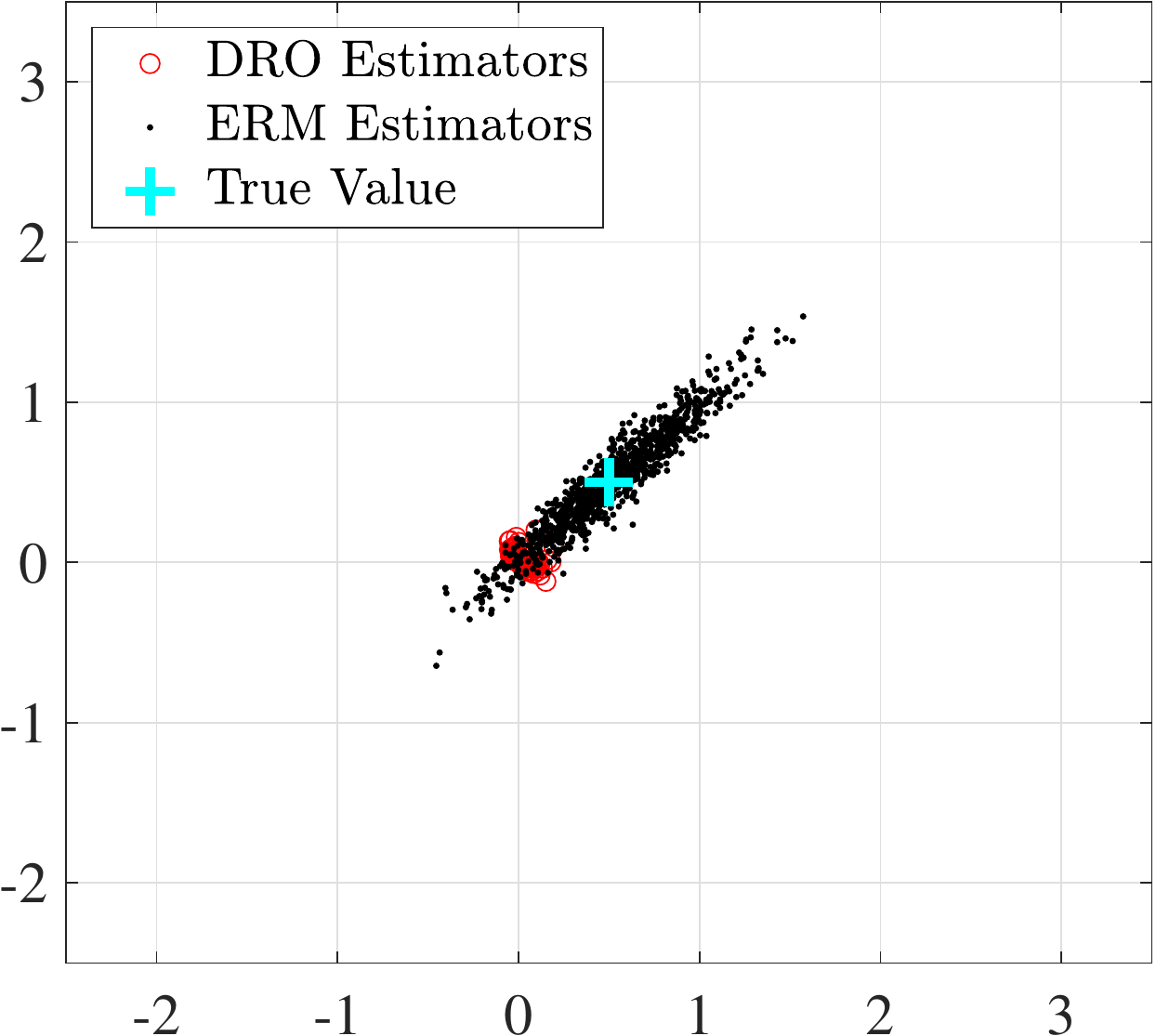}}
\caption{Scatter plots of $\ermin$ (black circles) and
  $\drmin$ (red circles) for $\theta_\ast = [0.5,0.5]^\T$.}
\label{fig:scatter-p5}
\end{figure}
\begin{figure}[ptbh]
\centering
\subfigure[$\rho=0.95$]{
 \includegraphics[width=1.5in]{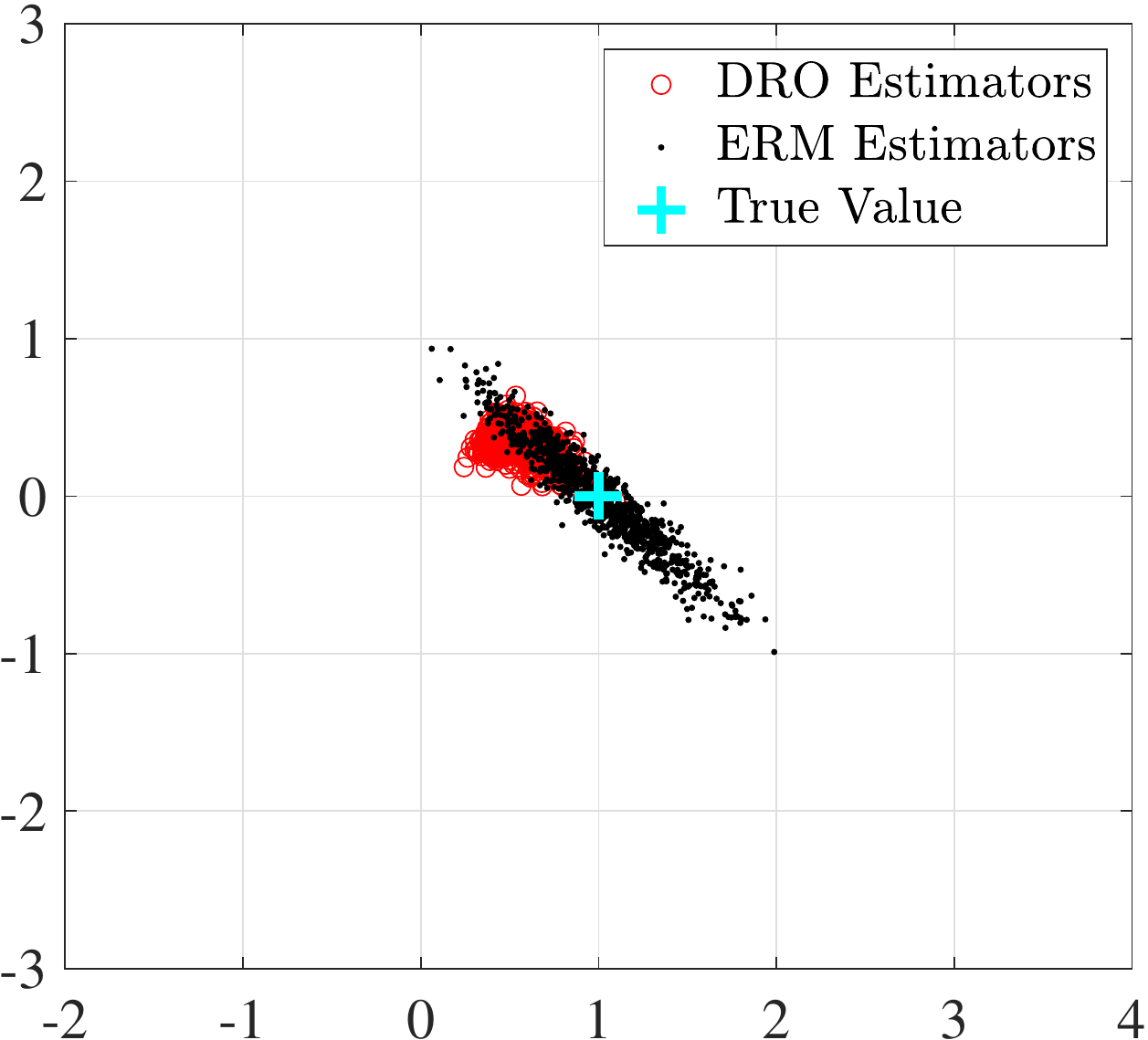}}
\subfigure[$\rho=0$]{
 \includegraphics[width=1.5in]{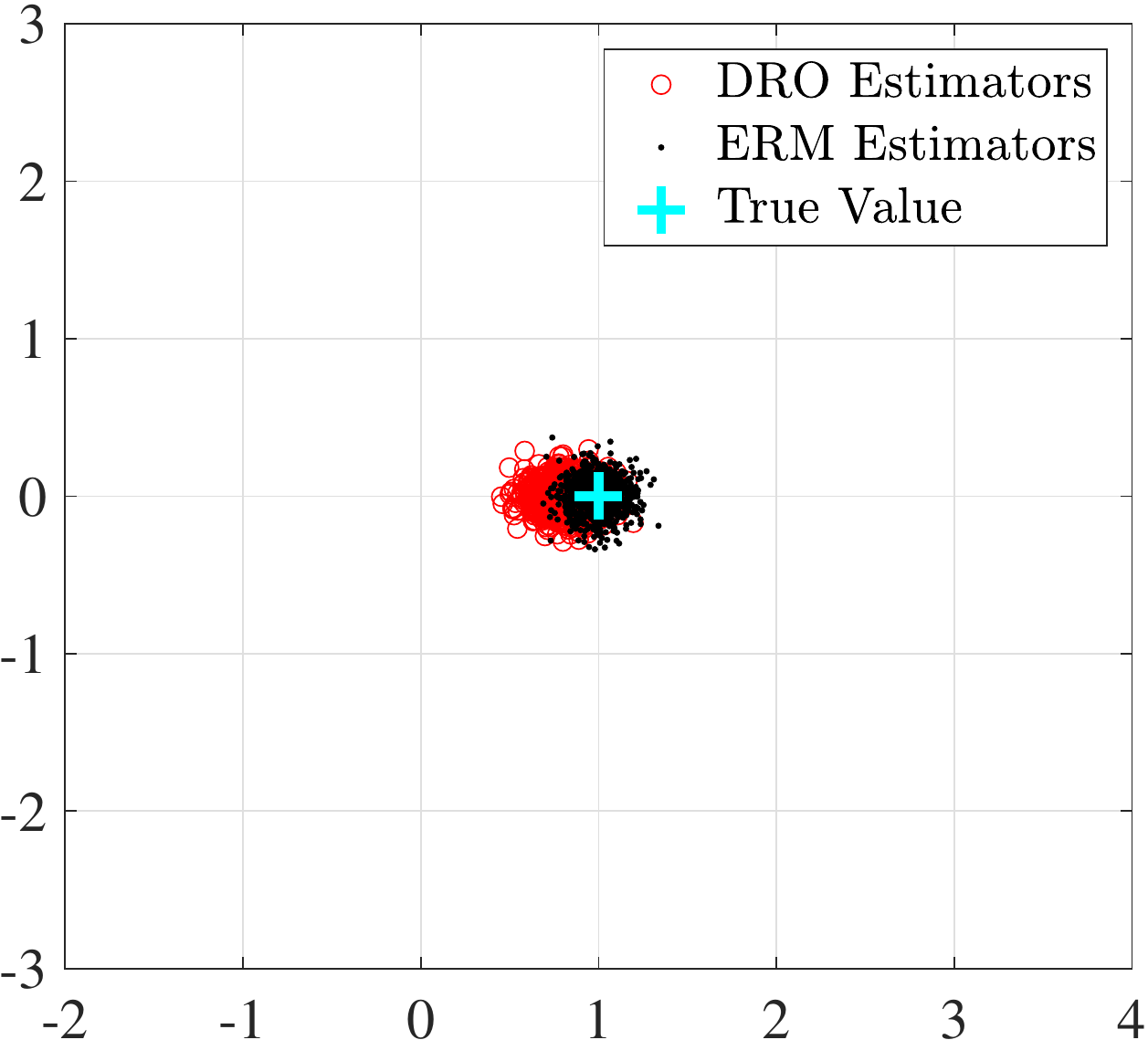}}
\subfigure[$\rho=-0.95$]{
\includegraphics[width=1.5in]{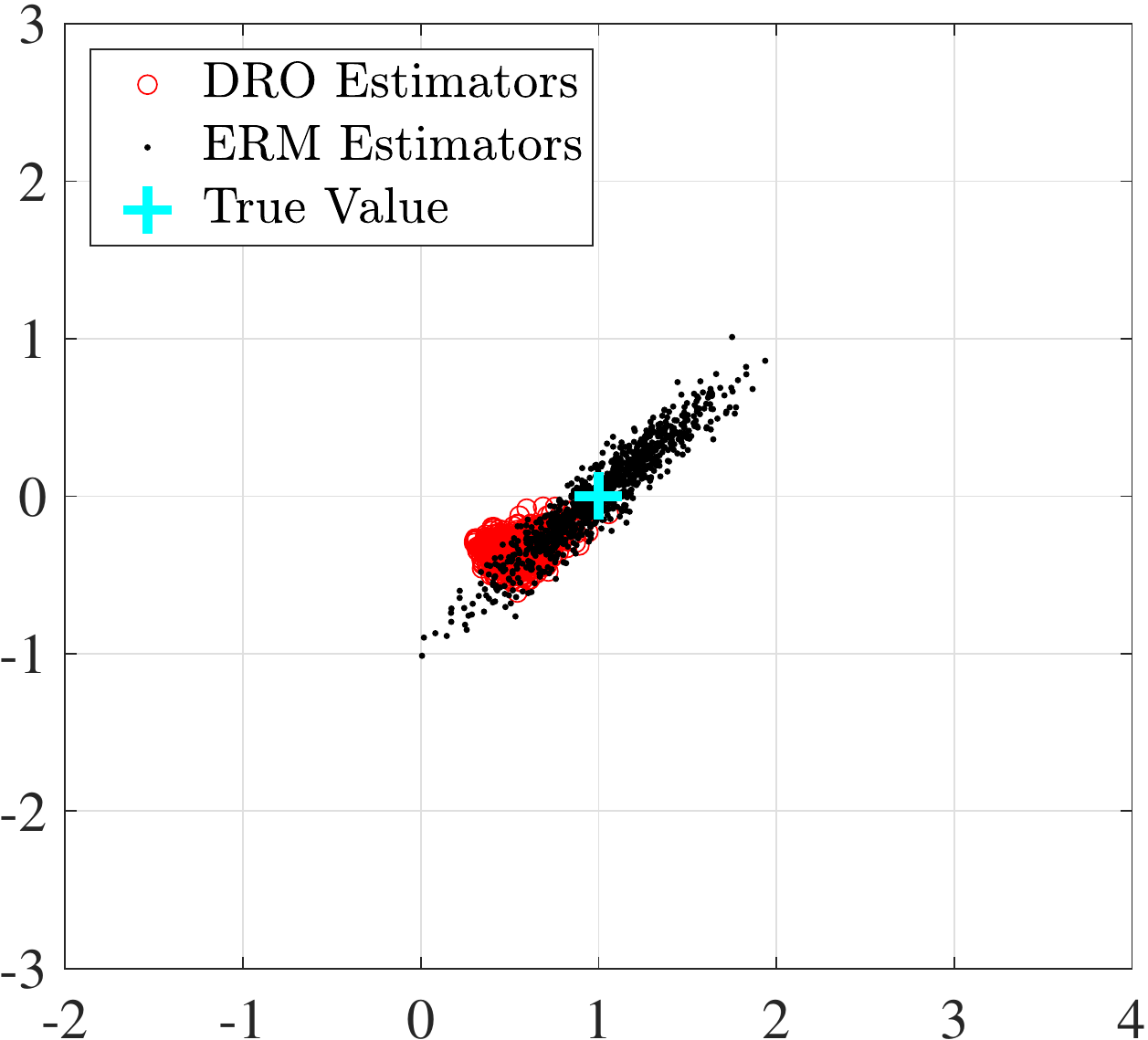}}
\caption{Scatter plots of $\ermin$ (black circles) and
  $\drmin$ (red circles) for $\theta_\ast = [1,0]^\T$. }
\label{fig:scatter-1}
\end{figure}
\label{eg:lin-reg}
\end{example}

\begin{remark}
\textnormal{Given the extensive coverage of the previous tutorial on Wasserstein DRO~\cite{ref:kuhn2019wasserstein} and of other extensive surveys on DRO in general such as~\cite{ref:rahimian2019distributionally}, we focus in this tutorial only on selected examples that we will use to connect to various statistical methodologies (including high-dimensional statistics and output analysis). However, DRO estimators have also been used in the context of non-parametric estimation (leading to the best known statistical rates) for instance, in the context of convex regression in which non-parametric functional regularization arises~\cite{ref:blanchet2019multivariate}. Other instances of regularization in the context of the 1-Wasserstein distance include support vector machines, norm-regularized
logistic and quantile regression~\cite{ref:shafieezadeh2017regularization}. Shrinkage behaviors of statistical estimators can also be induced by the 2-Wasserstein DRO estimators; see, for example,~\cite{ref:nguyen2018distributionally}. At the other end of the spectrum, the $\infty$-Wasserstein distance can be used to formulate robust conditional expectation or quantile estimators~\cite{ref:nguyen2020distributionally} and tackle multistage problems \cite{bertsimassturt}.} \label{rem:Lip-reg}
\end{remark}

\section{Informing the ambiguity radius $\delta$ via Optimal
Transport Projections}
\label{sec:WPF-Intro-Dev} 

We next consider the question of selecting the
ambiguity radius $\delta$ in \eqref{u-delta-defn} such that the resulting
DRO formulation has desirable statistical properties. A possible approach
towards this end is to select $\delta$ large enough so that the true data
generating process belongs to a distributional ambiguity set with some
prespecified confidence. This approach, which is often advocated in the
literature in machine learning and robust control~\cite{ref:shafieezadeh2015distributionally,ref:esfahani2018data, ref:obloj2021robust,ref:esfahani2018inverse,ref:yang1,ref:coulson2020distributionally,ref:8262716,ref:8910389,ref:9157949,ref:hakobyan2020wasserstein,boskos2021high,jiang2019data}, leads to a pessimistic
selection of $\delta$ simply because this criterion is not informed at all
by the loss function defining the decision problem. Indeed, the dimensional
dependence in the the concentration inequalities used for this purpose is
such that we will require an exponential amount in $m$ (the dimension of $X$) more samples to halve the error in the
resulting DRO solution (see the discussion after Theorem \ref{thm:joint-lim}
below).

Another approach involves the use of generalization bounds that are derived
to obtain finite sample performance guarantees. 
While appealing, this method usually requires distributions with compact support or sub-Gaussian tail assumptions on the underlying distributions and some of these bounds, in turn, rely on the convergence rate of the empirical Wasserstein distance~\cite{ref:fournier2015rate,  ref:trillos2015rate}. Another complicating factor in their use is that the multiplicative constants involved in the bounds tend to be pessimistic or difficult to compute for the
problem-in-hand.  

The most popular approach used in practice is based on cross-validation (CV). Despite its popularity, CV is often used in a way which could lead to inconsistent estimation (i.e., the incorrect identification of the optimal decision). For example, holdout and leave-one-out CV does not guarantee consistency in multivariate regression estimation. The $k$-fold CV approach leads to consistent estimation, but it requires $k/n \rightarrow 1$ and $n-k \rightarrow \infty$ as $n\rightarrow \infty$ (see, for example,~\cite{ref:shao1993linear}). As a consequence, when applying $k$-fold cross validation one needs to solve $k$ optimization problems. Typically, $k$ is chosen as a small number such as $k = 5$, but given the cautionary results in \cite{ref:shao1993linear} it is unclear if these choices are always appropriate relative to a given sample size $n$.

In the specific context of the distributionally robust precision matrix estimation \cite{ref:nguyen2018distributionally}, the Wasserstein ambiguity size that minimizes the \textit{expected} distance between the estimator and the true precision matrix scales linearly with the sample size as~$n^{-1}$, where the proportionality constant is a function of the true covariance matrix that is known in closed form \cite[Theorem~1]{ref:blanchet2019optimal}. The analysis in~\cite{ref:blanchet2019optimal}, however, depends on the analytical solution of the estimator and thus can only be generalized on a case-by-case basis.

Here we explain a generically applicable projection based statistical inference method called the 
\textit{Wasserstein Profile Function}, introduced in \cite{ref:blanchet2016robust}, for optimally informing the ambiguity radius $%
\delta$ in \eqref{u-delta-defn}. 
While the approach, as we shall see, enforces an optimality criterion at the level of decisions and not the value function, the methodology can be used to infer bounds on the optimal value as we shall discuss, see for example Theorem 4 in \cite{ref:blanchet2016robust}. Intuitively, the Wasserstein Profile Function computes the projection, along with the corresponding distance, of the empirical distribution $\Pr_n$ to a linear manifold of distributions which characterize the optimal decision. The use of these types of projection criteria to perform statistical inference is prevalent in statistics: \cite{ref:owen2001empirical} provides a comprehensive reference for projections, or profile functions, computed based on likelihood ratio metrics or the Kullback-Liebler divergence. The work of  \cite{ref:lam2017empirical,ref:duchi2021statistics,ref:lam2019recovering} connect these types of projections to DRO in the context of divergence measures.

Before giving the rationale behind the use of the Wasserstein Profile Function in the selection of $\delta$, we
first develop an understanding of this projection based inference procedure
in the elementary context of statistical hypothesis testing. The optimality
of the prescription of $\delta$ and related statistical implications are
explored in subsequent sections.

\subsection{Statistical hypothesis testing with projection based profile
inference.}
\label{sec:hyp-testing} Given independent samples $X_1,\ldots,X_n$ from the
unknown distribution $P_\ast$, we are interested in assessing if a given $\theta_0 \in \Theta$ satisfies the equation $E _{\Pr_\ast}[h(X,\theta_0)]
= 0$. Towards this end, for each $\theta \in \Theta$, let 
\begin{align}
\mathcal{F}_\theta\coloneqq \left\{ \Pr \in \mathcal{P}(\mathcal{X} %
): E _{\Pr}\left[ h(X,\theta) \right] = 0 \right\} \label{eq:F_theta}
\end{align}
denote the set of distributions of the random vector $X$ that satisfies the
condition $E_P[h(X,\theta)] = 0$. With this notation, our testing problem can be written as a 
statistical test with the hypotheses 
\begin{align*}
\mathcal{H}_0: \Pr_\ast \in \mathcal{F}_{\theta_0} \quad \text{ against } \quad 
\mathcal{H}_1: \Pr_\ast \notin \mathcal{F}_{\theta_0}.
\end{align*}
With the null hypothesis $\mathcal{H}_0$ stipulating the condition $E %
_{\Pr_\ast}[h(X,\theta_0)] = 0$, the statistical test will detect the failure
of the parameter choice $\theta_0$ in satisfying this condition with a pre-specified
confidence. Equipped with the Wasserstein distance, testing the
inclusion $P_\ast \in \mathcal{F}_{\theta_0}$ is equivalent to testing the
distance from $P_\ast$ to the set $\mathcal{F}_{\theta_0}$ being zero. With this
perspective, the hypotheses can be expressed as
\begin{equation*}
\mathcal{H}_0: \inf_{P \in \mathcal{F}_{\theta_0}}D_c(\Pr_\ast, P) = 0 \quad 
\text{ against } \quad \mathcal{H}_1: \inf_{P \in \mathcal{F}_{\theta_0}}
D_c(\Pr_\ast, P) > 0. 
\end{equation*}
To develop a suitable test statistic, we define the projection distance
function of the empirical distribution $\Pr_n = \frac{1}{n}\sum_{i=1}^n
\delta_{X_i}$ onto $\mathcal{F}_{\theta}$ as 
\begin{align}
\mathcal{P}(\Pr_n,\theta) \coloneqq \inf_{\Pr \in \mathcal{F}_\theta}
D_c(\Pr_n,\Pr) = \left\{ 
\begin{array}{cl}
\inf & D_c(\Pr_n,\Pr) \\ 
\mathrm{s.t.} & E _{\Pr}\left[ h(X,\theta)\right] = 0.%
\end{array}
\right. \label{eq:RWPI}
\end{align}
We address the projection metric $\mathcal{P}(P_n, \theta)$,
viewed as a function of $\theta$, as the \textit{Wasserstein Profile Function}. Equipped with the definition \eqref{eq:RWPI}, the statistical test will proceed generically as follows: for a
pre-specified significance level $\alpha \in (0, 1)$,

\begin{center}
reject $\mathcal{H}_0$ if $s_{n} > \eta_{_{1-\alpha}}$,
\end{center}
where $s_{n}$ is a test statistic that depends on the projection distance $%
\mathcal{P}(P_n, \theta_0)$, and $\eta_{_{1-\alpha}}$ is the $(1-\alpha)\times
100\%$ quantile of a limiting distribution obtained by studying the limit of 
$\mathcal{P}(P_n, \theta_0)$ as the number of samples $n$ tends to infinity.
To operationalize this statistical test, we will now examine a dual
reformulation of the projection distance $\mathcal{P}(P_n, \theta_0)$ and its
limiting behavior as $n \rightarrow \infty$. 
\subsubsection{Dual reformulation for $\mathcal{P}(P_n,\theta)$}

Similar to the dual formulation of the worst-case risk $\mathcal{R}_{\delta} %
(P_n,\theta)$ presented in Theorem \ref{thm:duality}, the projection
distance $\mathcal{P}(P_n,\theta)$ admits a dual reformulation due to the
linear programming structure offered by the optimal transport costs. To
state the dual reformulation, define $\tilde{\Theta} \subseteq \mathbb{R} %
^d$ as 
\begin{align*}
\tilde{\Theta} \coloneqq \left\{\theta \in \Theta^\circ: 0 \in \textnormal{%
conv}\left[ \{h(x,\theta) : x \in \mathcal{X} \}\right]^\circ\right\},
\end{align*}
where $\textnormal{conv}(S)$ denotes the convex hull of the set $S$, and $%
A^\circ$ indicates the interior of a set $A$. Since $\theta \notin \tilde{%
\Theta}$ cannot be a solution to $E _{\Pr}[h(X,\theta)] = 0$ unless $%
\Pr$ is degenerate, it is sufficient to restrict the analysis to $\theta
\in \tilde{\Theta}$.

\begin{proposition}[Dual reformulation for $\mathcal{P}(\Pr_n,\theta)
$]
\label{prop:dual-wpf} Let $h:\mathbb{R} ^m \times \mathbb{R} %
^d \rightarrow \mathbb{R} ^k$ be Borel measurable and $\{(x,x^\prime)
\in \mathcal{X}  \times \mathcal{X} : c(x,x^\prime) < +\infty
\}$ be Borel measurable and nonempty. Then for any $\theta \in \tilde{\Theta}%
$, 
\begin{align*}
\mathcal{P}(\Pr_n,\theta) = -\sup_{\lambda \in \mathbb{R} ^k} \frac{1%
}{n} \sum_{i=1}^n \sup_{x \in \mathcal{X} } \left\{ \lambda^\intercal %
 h(x,\theta) - c(X_i,x)\right\}.
\end{align*}
\end{proposition}

Verification of Proposition \ref{prop:dual-wpf}, as illustrated in \cite[Appendix B]{ref:blanchet2016robust}, utilizes semi-infinite linear program
duality. One may however quickly check the identity for the case where $\mathcal{X} $ is finite, as in the passage following Theorem \ref{thm:duality}.

\subsubsection{The limiting behavior of $\mathcal{P}(\Pr_n,\theta_0)$
under $\mathcal{H}_0$}

To see how we can obtain a test statistic from the projection $\mathcal{P}%
(\Pr_n,\theta_0)$, we next consider a suitably scaled version of $\mathcal{P}%
(\Pr_n,\theta_0)$ and obtain its limiting distribution when the
data-generating $\Pr_\ast \in \mathcal{F}_{\theta_0}$. Taking the special case $%
\mathcal{X}  = \mathbb{R} ^m$ and $c(x,x^\prime) = \Vert x -
x^\prime \Vert_q^2$, we rewrite the dual reformulation in Proposition \ref%
{prop:dual-wpf} as
\begin{align}
&n \times \mathcal{P}(\Pr_n,\theta_0)
\nonumber\\
&\qquad= \sup_{\lambda \in \mathbb{R} %
^k} \left\{ n^{1/2}\lambda^\intercal  H_n - \frac{1}{n} \sum_{i=1}^n
\sup_{x \in \mathbb{R} ^m} \left\{ n\lambda^\intercal  \big[ %
h(X_i + n^{-1/2} \Delta,\theta_0) - h(X_i,\theta_0)\big] - \Vert \Delta
\Vert_q^2 \right\} \right\}  \label{rwp-dual-rewrite} \\
&\qquad\approx \sup_{\xi \in \mathbb{R} ^k} \left\{ \xi^\intercal  %
H_n - \Phi_n(\xi,\theta_0)\right\},  \notag
\end{align}
where 
\begin{align*}
H_n \coloneqq n^{-1/2} \sum_{i=1}^nh(X_i,\theta_0) \quad \text{ and } \quad
\Phi_n(\xi,\theta) \Let \frac{1}{n}\sum_{i=1}^n \sup_{\Delta \in \mathbb{R} %
^m} \left\{ \xi^\intercal  D_x h(X_i, \theta) \Delta - \Vert
\Delta \Vert_q^2 \right\},
\end{align*}
and $D_x$ denotes the partial derivative with respect to the variable $x$.
The latter expression is obtained by changing variables from $\lambda$ to $%
\xi = \lambda n^{1/2}$, the variable $x$ inside the respective supremum as
in $x = X_i + n^{-1/2}\Delta$, and using a first-order approximation for the
difference $h(X_i + n^{-1/2} \Delta,\theta_0) - h(X_i,\theta_0)$. A precise
treatment of this term using the fundamental theorem of calculus can be found in 
\cite[Appendix A]{ref:blanchet2016robust}. 

The most important part of the derivation is the intuition provided by the scaling. The re-scaling involving $n^{-1/2}\Delta$ indicates that the optimal projection involves an optimal transport displacement of order $O(n^{-1/2})$, which is consistent with the Central Limit Theorem, but at a local level (meaning for each $X_i$).

Next, due to a similar application of H\"{o}lder's inequality as in Example \ref{eg:wc-return}, the suprema in the definition of $\Phi_n(\xi,\theta)$ can be simplified to result in $\Phi_n(\xi,\theta) = \frac{1}{4n} \sum_{i=1}^n \Vert \xi^\intercal  %
D_xh(X_i,\theta) \Vert_p^2$. For the sake of brevity, define 
\begin{align}
\varphi(\xi,\theta) \coloneqq \frac{1}{4} E _{\Pr_\ast} \Vert
\xi^\intercal  D_xh(X_i,\theta) \Vert_p^2 \quad \text{ and } \quad
\varphi^\ast(z,\theta) \coloneqq \sup_{\xi}~\{ \xi^\intercal  z -
\varphi(\xi,\theta)\}.  \label{phi-phist-defn}
\end{align}
Notice that $\varphi$ is a convex function in $\xi$, and $\varphi^*$ is the
convex conjugate of $\varphi$ with respect to the $\xi$ component. Then as a
consequence of the law of large numbers and uniform locally Lipschitz conditions we have 
\begin{align*}
\Phi_n(\xi,\theta) \rightarrow \varphi(\xi,\theta) \quad \text{ as } n
\rightarrow \infty
\end{align*}
uniformly over compact subsets of $\xi \in \mathbb{R} ^k$ and $\theta \in \Theta$. As a
result, 
\begin{align*}
n \times \mathcal{P}(P_n,\theta_0) = \varphi^\ast(H_n,\theta_0) (1 + o(1))
\end{align*}
as the number of samples $n \rightarrow \infty$. If $\theta_0$ indeed satisfies $E _{\Pr^\ast}[h(X,\theta_0)] = 0$, then the sequence $\{H_n\}_{n \geq 1}$ converges in distribution (due to the Central Limit
Theorem) and we obtain the limiting result as a consequence. In the
following, we use $\overset{dist.}{\longrightarrow} $ to denote
convergence in distribution.

\begin{theorem}[Limit theorem for $\mathcal{P}(\Pr_n,\theta)$]
\label{thm:limit-wpf} Suppose the function $h(\cdot,\theta_0)$ is continuously
differentiable and $E _{\Pr_\ast}[D_x h(X,\theta_0)D_x
h(X,\theta_0)^\intercal ] \succ 0$. Then under the null hypothesis $%
\mathcal{H}_0$, 
\begin{align*}
n \times \mathcal{P}(P_n,\theta_0) \overset{dist.}{\longrightarrow}  %
\varphi^\ast(H,\theta_0) \quad \text{ as } n \rightarrow \infty,
\end{align*}
where $H \sim \mathcal{N}(0,\mathrm{Cov}_{P_\ast}[h(X,\theta_0)])$.
\end{theorem}

\begin{remark}
\textnormal{Recall that $c(x,x^\prime) = \Vert x -
x^\prime \Vert_q^2$, so in terms of the Wasserstein distance of order 2 (i.e., $r=2$) the projection converges to zero under the null hypothesis at rate $O(n^{-1/2})$. The proof of this result, along with extensions for the case $r \geq 1$, is given in \cite{ref:blanchet2016robust} assuming that $(X_1,...,X_n)$ is an i.i.d.~collection. As one can see from the discussion leading to Theorem \ref{thm:limit-wpf}, the key ingredients really are a functional law of large numbers for $\mathcal{P}(P_n,\theta)$ over compact sets and a Central Limit Theorem for $H_n$, both of which hold well beyond the i.i.d.~assumptions imposed here and in \cite{ref:blanchet2016robust}. The work \cite{ref:blanchet2020distributionally} discusses non-i.i.d.~extensions which are relevant in financial applications, in particular, the mean-variance portfolio allocation problem.}
\end{remark}
Conceptually, Theorem \ref{thm:limit-wpf} reveals that $s_n = n\times%
\mathcal{P}(\Pr_n,\theta_0)$ serves as a test statistic to reject the null
hypothesis $\mathcal{H}_0$. In particular, for a pre-specified significance
level $\alpha \in (0, 1)$, let $\eta_{_{1-\alpha}}$ denote the $(1 - \alpha)
\times 100\%$ quantile of the limiting distribution given by the law of $%
\varphi^\ast(H,\theta_0)$. Then rejecting $\mathcal{H}_0$ if $s_n >
\eta_{_{1-\alpha}}$ results in a statistical test with Type-I error probability $\alpha$. The test distribution is determined by $\text{Cov}_{\Pr_\ast}[h(X,\theta_0)]$. The test distribution is unaffected by using any consistent plug-in estimator of the covariance matrix if it is unknown.

We now dive into an application of the proposed hypothesis test for fair classification~\cite{ref:taskesen2020statistical}.

\begin{example}[Test for probabilistic fair classifier]
\textnormal{
    Consider the joint random vector $(X, A, Y)$ consisting of a feature vector $X \in \R^d$, a sensitive attribute $A \in \{0, 1\}$ and a class label $Y \in \{0, 1\}$. A logistic classifier aims to predict the class label $Y$ for any feature input $X$. A logistic classifier can be represented by the conditional distribution function $h_\theta$ of $Y$ given $X$ of the sigmoid form
    \[
        h_\theta(x) = \frac{1}{1 + \exp(-\theta^\top x)}.
    \]
    Following the definition in~\cite{ref:pleiss2017fairness}, we say that a probabilistic classifier $h_\theta$ satisfies the probabilistic equal opportunity criterion relative to a distribution $\Pr$ if
    \[
        E_{\Pr}[h_\theta(X) | A = 1, Y = 1] = E_{\Pr}[h_\theta(X) | A = 0, Y = 1].
    \]
    As a consequence, the manifold of distributions that renders $h_\theta$ a fair classifier is defined specifically as
    \[
     \mathcal F_{\theta}= \left\{ 
     Q \in \mathcal{P}(\R^d \times \{0, 1\} \times \{0, 1\}):
     E_{Q}[h_\theta(X) | A = 1, Y = 1] = E_{Q}[h_\theta(X) | A = 0, Y = 1]
    \right\}.
    \]
    Suppose that $(X_i, A_i, Y_i)$ are i.i.d.~samples from $P_\ast$ and that the ground transport cost is of the form $ c\big( (x', a', y'),  (x, a, y) \big) = \| x - x'\|_q^2 + \infty | a - a'| + \infty | y - y'|$. Under the null hypothesis $\mathcal H_0:~P_\ast \in \mathcal{F}_\theta$, we have the following limit distribution
    \begin{align*}
        &n \times \mathcal{P}(\Pr_n, \theta)  \dconv \beta \chi_1^2,
    \end{align*}
    where $\chi_1^2$ is a chi-square distribution with 1 degree of freedom,
    \[
    \beta = \left(E_{P_\ast} \left[\left\| \nabla h_\theta(X) \left( \frac{\mathbb{I}_{(1, 1)}(A, Y)}{ p_{11}}- \frac{\mathbb{I}_{(0, 1)}(A, Y)}{p_{01}} \right) \right \|_p^2 \right] \right)^{-1} \frac{\sigma^2}{p_{01}^2 p_{11}^2}
    \] 
    with $\sigma^2 = \mathrm{Cov}( Z)$, $p_{a1} = \Pr_\ast(A = a, Y = 1)$, $\mathbb{I}$ is the indicator function and $Z$ is the univariate random variable
    \begin{align*}
        Z&= h_\theta(X)\left( p_{01}\mathbb{I}_{(1,1)}( A,Y) 
     -p_{11}\mathbb{I}_{(0,1)}( A, Y) \right)\\ & \qquad +\mathbb{I}_{(0,1)}(
    A,Y) E_{P_\ast}[\mathbb{I}_{(1,1)}( A, Y)
    h_\theta(X)] -\mathbb{I}_{(1,1)}( A, Y) E_{P_\ast}[\mathbb{I}_{(0,1)}(
    A, Y) h_\theta(X)].
\end{align*}
Given a logistic classifier parametrized by $\theta$, the decision to reject the probabilistic fairness of this classifier now relies on computing the projection distance $\mathcal P(P_n, \theta)$ and on computing the empirical quantile estimate of $\beta \chi_1^2$~\cite{ref:taskesen2020statistical}. 
} \label{eq:fair}
\end{example}

The limit result in Theorem~\ref{thm:limit-wpf} relies on the assumption that $h$ is continuously differentiable. The extension of the limiting distribution when $h$ is non-differentiable, or even when $h$ is discontinuous, can be found in~\cite{ref:si2021testing}. The machinery of the Wasserstein Profile Function is also applied for model selection of graphical Lasso~\cite{ref:cisneros2020distributionally}. 
It is important to notice that the hypothesis testing framework outlined in this section is fundamentally different from the robust hypothesis test with Wasserstein ambiguity set proposed in~\cite{ref:gao2018robust}. Therein, the test is constructed to minimize the worst-case error, measured by the maximum of the type-I and type-II errors, uniformly over all perturbations of the empirical distribution in the ambiguity set.

The next section explains how this hypothesis testing procedure can guide us
to choose the ambiguity radius $\delta$ optimally in a certain statistical
sense.

\subsection{Informing DRO ambiguity radius $\delta$ from the
projection $\mathcal{P}(P_n,\theta)$}
Going back to the DRO formulation \eqref{dro}, assume $\ell(x,\cdot)$ is
convex for every $x \in \mathcal{X} $ and let
\begin{equation*}
h(x,\theta) \coloneqq D_\theta \ell(x,\theta)
\end{equation*}
denote the partial derivative of the loss function $\ell$. Fix an arbitrary
distribution $\Pr$ in the ambiguity set $\mathcal{U}_\delta(\Pr_n)$ and
suppose that $E _{\Pr}[h(X,\theta)] = 0$ specifies the necessary and
sufficient condition for minimizing the risk $\mathcal{R} (\Pr,\theta)
$ over the feasible parameters $\theta \in \Theta$. In this case, the set $%
\{\theta \in \Theta : \Pr \in \mathcal{F}_\theta\}$ (recall the definition of $\mathcal{F}_\theta$ in \eqref{eq:F_theta}) contains all parameter
choices that are optimal from the decision maker's point of view.
Consequently, by taking unions over all $P \in \mathcal{U}_\delta(\Pr_n)$,
the set 
\begin{align}
\Lambda_\delta(\Pr_n) &\coloneqq \left\{ \theta \in \Theta: \mathcal{F}%
_\theta \cap \,\mathcal{U}_\delta(\Pr_n) \neq \varnothing \right\}
\label{Lambda-defn} \\
&= \left\{ \theta \in \Theta: \mathcal{R}  (\Pr,\theta) = \mathcal{R} %
_\ast(\Pr)\right\}  \notag
\end{align}
includes all the parameter choices that are collected by the decision maker
as optimal for some distribution in the distributional uncertainty set $%
\mathcal{U}_\delta(\Pr_n)$. This leads to the following notion of compatible confidence regions. 

\subsubsection{Confidence regions compatible with the DRO formulation \eqref{dro}.}
\label{sec:intro-cc}
If $\mathcal{U}_\delta(\Pr_n)$ represents a family of \textit{plausible} representations of uncertainty around $P_n$, then $\Lambda_\delta(\Pr_n)$ in \eqref{Lambda-defn} represents a family of plausible decisions. One can therefore think of $\Lambda_\delta(\Pr_n)$ as the projection of $
\mathcal{U}_\delta(\Pr_n)$ onto the decision space (which is finite-dimensional and informed by the optimization problem of interest). In that sense, one can view $\Lambda_\delta(\Pr_n)$ as a set of decisions $\theta \in \Theta$ that are \textit{compatible} with the distributional uncertainty $\mathcal{U}_\delta(\Pr_n).$ If we can guarantee that the set $\Lambda_\delta(\Pr_n)$ contains an optimal $\theta_\ast$ solving \eqref{prm} with probability $1-\alpha,$ then $\Lambda_\delta(\Pr_n)$ becomes a compatible confidence region for a correct decision for the problem \eqref{prm}.

Since the family of sets $\{\Lambda_\delta(\Pr_n) : \delta > 0\}$ is increasing in $\delta$ it is clear that a minimizer of \eqref{prm} will be a member of $%
\Lambda_\delta(\Pr_n)$ for a sufficiently large choice of $\delta$. Indeed,
this holds true when $\Pr_\ast \in \mathcal{U}_\delta(\Pr_n)$, or more
explicitly, when $\delta$ is chosen such that $D_c(\Pr_n, \Pr_\ast) \le
\delta$. Based on these considerations, since solving \eqref{prm} is our primary objective, this leads to
the following natural enquiry:\newline

\begin{tcolorbox}[colback=white!5!white,colframe=black!75!black]
  \textbf{Question:} What is the smallest $\delta > 0$ such that the confidence region $\Lambda_\delta(\Pr_n)$ contains an unknown population risk minimizer of \eqref{prm} with a target $(1-\alpha)$-confidence? 
\end{tcolorbox}
In other words, we seek to identify 
\begin{align}
\delta_\ast \coloneqq \inf\big\{ \delta > 0: \Lambda_\delta(\Pr_n) \text{
contains a minimizer of \eqref{prm} with } (1 - \alpha)\text{-confidence} %
\big\}.   \label{delta-st-defn}
\end{align}
The desire to identify the smallest $\delta$ (satisfying this criterion) is
motivated by the need to drive down the conservativeness of the resulting
DRO formulation \eqref{dro}.

\subsubsection{Identifying the optimal radius $\delta_\ast$ from the projection profile $\mathcal{P}(\Pr_n,\theta).$}
The above question leads to a data-driven choice of the ambiguity size $\delta$
that is explicitly linked to the statistician's decision problem and can be
answered with the projection-based profile function $\mathcal{\Pr}%
(\Pr_n,\theta)$ introduced in Section \ref{sec:hyp-testing} for testing the
hypothesis $\Pr_\ast \in \mathcal{F}_\theta$. To see the link between these
apparently different exercises, suppose that any solution $\theta_\ast \in 
\text{arg}\,\text{min} _{\theta \in \Theta} E _{\Pr_\ast}\left[
\ell(X,\theta)\right]$ satisfies the optimality condition
\begin{align*}
E _{\Pr_\ast} \left[ h(X,\theta_\ast) \right] = 0. 
\end{align*}
Then the minimal radius $\delta_\ast$ in~\eqref{delta-st-defn} can be
re-expressed as 
\begin{align*}
\delta_\ast &= \inf\left\{\delta > 0: \exists \theta \in
\Lambda_\delta(\Pr_n) \text{ such that } \Pr_\ast \in \mathcal{F}_\theta 
\text{ with } (1 - \alpha) \text{-confidence}\right\} \\
&\overset{(a)}{=} \inf \left\{ \delta > 0: \exists (\Pr,\theta) \in \mathcal{%
U}_\delta(\Pr_n) \times \Theta \text{ such that } \Pr \in \mathcal{F}%
_\theta, \Pr_\ast \in \mathcal{F}_\theta \text{ with } (1 - \alpha) \text{%
-confidence} \right\} \\
&\overset{(b)}{=} \inf\left\{\delta > 0: \exists \theta \in \Theta \text{
such that } \mathcal{P}(\Pr_n,\theta) \leq \delta, \ \Pr_\ast \in \mathcal{F}%
_\theta \text{ with } (1 - \alpha) \text{-confidence} \right\} \\
&\overset{(c)}{=} \inf \left\{\delta > 0: \exists \theta \in \Theta \text{
such that } \mathcal{P}(\Pr_n,\theta) \leq \delta, \ n\mathcal{P}%
(\Pr_n,\theta) \leq \eta_{_{1-\alpha}}(1+o(1)) \right\} \\
&= n^{-1} \times \eta_{_{1-\alpha}}(1+o(1)),
\end{align*}
where (a) and (b) follow respectively from the definitions of $%
\Lambda_{\delta}(\Pr_n)$ and of $\mathcal{P}(\Pr_n,\theta)$, (c) from the
statistical test developed in Section \ref{sec:hyp-testing} for rejecting
the null hypothesis $\mathcal{H}_0: \Pr_\ast \in \mathcal{F}_\theta$ with $%
(1-\alpha)$-confidence, and the last equality holds if a solution $%
\theta_\ast$ to \eqref{prm} satisfying $\Pr_\ast \in \mathcal{F}%
_{\theta_\ast}$ exists.

Thus, instead of requiring $\delta > 0$ to be large enough such that the
data-generating $\Pr_\ast \in \mathcal{U}_\delta(\Pr_n)$, the projection
based prescription merely requires existence of $\Pr \in \mathcal{U}_\delta$
such that $\Pr \in \mathcal{F}_{\theta}$ for some population risk minimizer $%
\theta_\ast \in \Theta$. While this choice results in $\delta \propto n^{-1}$
and the guarantee that $\Lambda_\delta(\Pr_n)$ serves as a confidence region
for the solutions to \eqref{prm}, the former prescription from concentration
inequalities results in a pessimistic $\delta \propto n^{-2/m}$ under
additionally restrictive assumptions, where $m$ is the ambient dimension of the space $\xsp$ in which the random vector $X$ takes values (see, for example, \cite{ref:weed2019sharp, ref:barrio1999central, ref:dudley1969speed, ref:kuhn2019wasserstein}).

Algorithm \ref{algo:recipe-dro} below provides a recipe for estimating $\delta$ in order to guarantee asymptotic optimality in the sense of ensuring the smallest choice which guarantees $1-\alpha$ coverage. 

\begin{algorithm}[h!]
  \caption{A recipe for DRO estimation \eqref{dro} with optimal selection of the ambiguity radius $\delta$.}   
  \begin{algorithmic}
  \STATE{\textbf{Input:} Training samples $X_1,\ldots,X_n$ from an unknown distribution $\Pr_\ast$. }\\
  \STATE{\textbf{Outline of steps:} } \\
  \STATE{1. Let $\hat{\Sigma}$ be the sample covariance of $\{h(X_1,\ermin),\ldots,h(X_n,\ermin)\}$, where $h(X,\theta) = D_\theta \ell(X,\theta)$ and $\ermin \in \argmin_{\theta \in \Theta} ~\E_{\Pr_n}[\ell(X,\theta)]$. }\\
   \STATE{2. Obtain independent samples  $H_1,\ldots,H_k$  from $\mathcal{N}(0,\hat{\Sigma})$. (Any choice of $k$ such that $k \rightarrow \infty$ as $n \rightarrow \infty$ is valid, for example, $k= \log(n)$.) Let $\hat{\eta}_{1-\alpha}$ be the $(1-\alpha)$-quantile of the sample collection $\{\hat{\varphi}^\ast(H_i,\ermin): i = 1,\ldots,k\}$, where $\hat{\varphi}^\ast$ is the convex conjugate of $\hat{\varphi}(\cdot)$ in \eqref{phi-savg}. In the absence of the knowledge of the conjugate $\hat{\varphi}^\ast$, one may take $\hat{\eta}_{1-\alpha}$ of a tractable upper bound of $\varphi^\ast(\cdot)$, as illustrated in Example \ref{eg:high-dim-lasso}. }\\
   \STATE{3. Let $\drmin \in \argmin_{\theta \in \Theta}~\drorisk(\Pr_n,\theta)$ be a solution to the DRO estimation \eqref{dro} obtained by letting $\delta = n^{-1} \times \eta_{1-\alpha}$.} \\
   \STATE{4. Return $\drmin$ as the distributionally robust solution.}
   \end{algorithmic}
    \label{algo:recipe-dro}
    \end{algorithm}

\subsubsection{Illustrative examples.} We illustrate the choice of $\delta$ prescribed in Algorithm \ref{algo:recipe-dro} via some illustrative examples. 
\begin{example}[Identifying the radius $\delta = n^{-1} \eta_{1-\alpha}$ for  linear regression example] \label{eg:lin-reg-delta}
Continuing the discussion in Example \ref{eg:lin-reg}, we see that $%
h(x,y,\theta) = D_\theta \ell(x,y,\theta) = (y - \theta^\intercal  x)x
$. Consider the linear regression model 
\begin{align}
Y = \theta_\ast^{\intercal } X + e,  \label{lin-reg-model}
\end{align}
where the unknown $\theta_\ast \in \mathbb{R} ^d$ and the independent
additive error $e$ has zero mean and variance $\sigma^2$. Letting $\Pr_\ast$
be the distribution of $(X,Y)$, we have $E _{\Pr_\ast}[h(X,Y,\theta_%
\ast)] = E _{\Pr_\ast} [e X] = 0$. Then under the null hypothesis $\mathcal{H}_0: \Pr_\ast \in \mathcal{F}_{\theta_\ast}$, we have from Theorem %
\ref{thm:limit-wpf} that $n\mathcal{P}(\Pr_n,\theta_\ast) \overset{dist.}{%
\longrightarrow}  \varphi^\ast(H,\theta_\ast)$, where $H \sim 
\mathcal{N}(0,\sigma^2 \Xi)$, $\Xi = E_{P_\ast}[XX^\intercal ]$ and 
\begin{align*}
\varphi^\ast(z,\theta_\ast) &= \sup_{\xi} \left\{ \xi^\intercal  z -
4^{-1} E _{\Pr_\ast} \left\Vert e \xi - \xi^\intercal  X
\theta_\ast \right\Vert_p^2 \right\}.
\end{align*}
Taking $q = 2$ in the transportation cost in \eqref{trans-cost-linreg} for
ease of illustration, we see that 
\begin{align*}
\varphi^\ast(z,\theta_\ast) &= \sup_{\xi} \left\{ \xi^\intercal  z -
4^{-1} \xi^\intercal  \left( \sigma^2 \mathbb{I}_d + \Xi \Vert
\theta_\ast \Vert_2^2 \right)\xi\right\} \\
&= z^\intercal  \left( \sigma^2 \mathbb{I}_d + \Xi \Vert \theta_\ast
\Vert_2^2 \right)^{-1} z.
\end{align*}
The limit $\varphi^\ast(H,\theta_\ast)$ has a generalized chi-square distribution in this case. 
If $\E_{\Pr_\ast}[X] =0$  and $\Xi$ is invertible with eigen decomposition $\Xi = UDU^\T$, then we have that $N = D^{-1/2}U^\T H$ is a standard normal vector with mean 0 and covariance $\mathbb{I}_d$. As a result, 
\begin{align*}
    \varphi^\ast(H,\theta_\ast) = \sum_{i=1}^d\frac{D_{ii}}{1+D_{ii}\Vert \theta_\ast\Vert_2^2/\sigma^2} N_i^2,
\end{align*}
where $D_{ii}$ is the $i$-th diagonal element of the matrix $D$. One may compute an estimate of the quantile $\hat{\eta}_{_{1-\alpha}}$ by plugging in any consistent estimators for $\theta_\ast$ and $\Xi$. The asymptotic validity of $\delta$ in terms of coverage and hypothesis testing remains unchanged by plugging in consistent estimators in this case because the limiting distribution function is continuous as a function of these parameters. 
\end{example}

\begin{example}[A prescription for $\delta$ in high-dimensional
linear regression] \label{eg:high-dim-lasso}
We now consider the linear regression model~\eqref{lin-reg-model} in the high-dimensional setting where $d \gg n$.
Choosing the transportation cost $c$ in \eqref{trans-cost-linreg} with $q =
a= \infty$, the resulting DRO linear regression problem coincides with the
square-root Lasso estimator in~\cite{ref:belloni2010square}. We now examine
the prescription $\delta = n^{-1} \eta_{_{1-\alpha}}$ arising from the
projection metric $\mathcal{P}(\Pr_n,\theta)$. Since $d \rightarrow \infty$
as $n \rightarrow \infty$, the limiting characterization in Theorem \ref%
{thm:limit-wpf} does not hold and therefore we take $\eta_{_{1-\alpha}}$
directly to be the $(1-\alpha)$-quantile of the pre-limit $n\mathcal{P}%
(\Pr_n,\theta_\ast)$. Resorting to \eqref{rwp-dual-rewrite} for this
purpose, we see that $n\mathcal{P}(\Pr_n,\theta_\ast)$ is merely the convex
conjugate of 
\begin{align*}
\tilde{\Phi}_n(\xi)&= \frac{1}{n} \sum_{i=1}^n \sup_{\Delta} \left\{ n^{1/2}
\xi^\intercal  \lbrack Y_i - \theta_\ast^\intercal  (X_i +
n^{-1/2} \Delta)(X_i + n^{-1/2}\Delta)] - \Vert \Delta \Vert_\infty^2
\right\} \\
&= \frac{1}{n} \sum_{i=1}^n \sup_{\Delta} \left\{ e_i \xi^\intercal  %
\Delta - (\theta_\ast^\intercal  \Delta)(\xi^\intercal  X_i) - %
\left[ 1 + n^{-1/2}\frac{(\theta_\ast^\intercal  \Delta)(\xi^%
\intercal  \Delta)}{\Vert \Delta \Vert_\infty^2} \right] \Vert \Delta
\Vert_\infty^2 \right\},
\end{align*}
evaluated at $H_n = n^{-1/2}\sum_{i=1}^n e_i X_i$. Bounding the
inner-products involving $\Delta$ using H\"{o}lder's inequality results in
the following nonasymptotic bound for the inner suprema: 
\begin{align*}
\tilde{\Phi}_n(\xi) \geq \Phi_n(\xi) \coloneqq \frac{1}{n} \sum_{i=1}^n 
\frac{\Vert e_i \xi - (\xi^\intercal  X_i) \theta_\ast\Vert_1^2}{1 +
n^{-1/2} \Vert \theta_\ast\Vert_1 \Vert \xi \Vert_1}.
\end{align*}
The convex conjugate of $\tilde{\Phi}_n(\cdot)$, denoted by $\Phi_n^\ast
(\cdot)$, can be bounded as $\Phi_n^\ast(z) \leq \Vert z \Vert_\infty^2/%
\text{Var}_n \vert e \vert$ $( 1 + o(1))$. Here $\text{Var}_n \vert e \vert$
is the sample variance of the collection $\{\vert e_i \vert\}_{i=1}^n$. The
intermediate steps involved in arriving at the upper bound of the convex
conjugate is available in the proof of Theorem 7 in \cite%
{ref:blanchet2016robust}. Since $n\mathcal{P}(\Pr_n,\theta_\ast) =
\sup_{\xi} \{ \xi^\intercal  H_n - \tilde{\Phi}_n(\xi)\}$, the
following upper bound for $n\mathcal{P}(\Pr_n,\theta_\ast)$ can be obtained
from an upper bound for the convex conjugate of $\Phi_n(\cdot):$ 
\begin{align*}
n\mathcal{P}(\Pr_n,\theta_\ast) \leq \frac{\Vert H_n \Vert_\infty^2}{\text{%
Var}_n \vert e \vert} (1 + o(1)),
\end{align*}
as $n \rightarrow \infty$. Suppose that the additive noise $e$ is normally
distributed and the observations $X_i = (X_{i1},...,X_{id})$ are normalized
so that $n^{-1}\sum_{i=1}^{n} X_{ij}^2 = 1$, for $j = 1,\ldots,d$. Then for
any $\alpha < 1/8, C > 0, \varepsilon > 0$, one can conclude the following
from \cite[Lemma 1(iii)]{ref:belloni2010square}: Conditional on the
observations $\{X_i : i = 1, . . . , n\}$, 
\begin{align*}
n\times \mathcal{P}(\Pr_n, \theta_\ast) \leq \eta_{_{1-\alpha}} \coloneqq %
\left[\frac{\pi}{\pi-2} \Phi^{-1} (1-\alpha/2d) \right]^2,
\end{align*}
with probability larger than $1-\alpha$, as $n \rightarrow \infty$,
uniformly in $d$ such that $\log d \leq C n^{1/2-\varepsilon}$. Here $%
\Phi^{-1}(1-\alpha)$ denotes the $(1-\alpha)$-quantile of the standard
normal random variable. This results in the choice of ambiguity radius $%
\delta = n^{-1} \times \eta_{_{1-\alpha}}$. The respective regularization
parameter in \eqref{sqroot-lasso} is given by
\begin{align*}
\delta_n^{1/2} = n^{-1/2} \times \frac{\pi}{\pi-2} \Phi^{-1} (1-\alpha/2d),
\end{align*}
which agrees with the prescription obtained independently in the statistics
literature for recovering $\theta_\ast$ in high-dimensional settings when $d
\gg n$; see, for example, \cite[Corollary 1]{ref:belloni2010square}. Another
interesting aspect of this choice is its self-normalizing property that
renders the selection independent of the error variance $\sigma^2$.
\end{example}

\section{Statistical properties of DRO estimators and the optimality of $%
\delta = n^{-1} \eta_{_{1-\alpha}}$}
\label{sec:dro-clt} 

For ease of illustrating the key ideas, we make the
following simplifying assumptions throughout this section. A precise statement of results in more general settings can be found in the accompanying references.

\begin{assumption}
The transportation cost  $c(x,x^\prime) = \Vert x - x^\prime \Vert_q^2$ for
some $q \in (1,\infty]$.  \label{assume:cost}
\end{assumption}

\begin{assumption}
The loss $\ell:\mathbb{R} ^m \times \Theta \rightarrow \mathbb{R} %
$ is twice continuously differentiable with uniformly bounded second
derivatives. For any $\theta \in \Theta$, $\{\ell(X,\theta): \theta \in
\Theta\}$ has finite second moments. \label{assume:loss}
\end{assumption}

\subsection{Adaptive regularization induced by Wasserstein DRO formulations}
\label{sec:adapt-reg} 

We begin our analysis of the DRO estimation problem~\eqref{dro}
with the series expansion of the DRO objective $\mathcal{R}_{\delta} %
(\Pr_n, \theta)$ in Theorem \ref{thm:var-reg} below. For $p$, $q$ satisfying $1/p
+ 1/q = 1$, let 
\begin{align}
\mathcal{V}^2 _p(\Pr,\theta) \coloneqq E _{\Pr}\Vert D_x\ell
(X,\theta)\Vert_p^2  \label{sens-defn}
\end{align}
denote the expected `squared variation' in the loss or, in other words, the expected squared sensitivity of loss with respect to perturbations in the random vector $X$. Theorem \ref{thm:var-reg}
below asserts that the DRO estimation procedure \eqref{dro} favours
solutions which possess low sensitivity to perturbations, measured in terms
of the regularizer $\mathcal{V} (\Pr_n,\theta)$. 

\begin{theorem}[Variation regularization]
Under Assumptions \ref{assume:cost} - \ref{assume:loss}, we have
\begin{align*}
\mathcal{R}_{\delta} (\Pr_n,\theta) = E _{\Pr_n} \left[
\ell(X,\theta)\right] + \delta^{1/2} \mathcal{V} _p(\Pr_n, \theta) +
O_p(\delta),
\end{align*}
for any choice of $\delta \rightarrow 0$ as $n \rightarrow \infty$. The
convergence is uniform over $\theta$ in compact subsets of $\Theta$. \label%
{thm:var-reg}
\end{theorem}

A precise characterization of the second-order error term, introduced by
carefully incorporating the second-order terms in the Taylor expansion of $%
\ell$, is available in \cite[Appendix A1]{ref:blanchet2021confidence}. A
general version of the result, applicable for $c(x,x^\prime) = \Vert x -
x^\prime \Vert_q^r$ (with $r$ is not necessarily required to equal 2), is
available in \cite{gaoChenKleywegt}. The limiting analysis in \cite{ref:blanchet2021confidence,ref:bartl2020robust} also characterizes the sensitivities of the optimizer as the radius of the Wasserstein ball is shrunk to zero. 

For small values of the ambiguity radius $\delta$, Theorem \ref{thm:var-reg}
asserts that the DRO objective $\mathcal{R}_{\delta} (\Pr_n,\theta)$
can be understood in terms of the empirical risk $\mathcal{R} %
(\Pr_n,\theta)$ and a regularization term $\mathcal{V} (\Pr_n,\theta)$
capturing the variation/sensitivity induced by the choice $\theta$. This
connection with the empirical risk minimization objective allows the
following interpretation: If there are several solutions with small
empirical risk (which happens often in high-dimensional settings), the DRO
estimation procedure \eqref{dro} can be understood as favouring solutions
which possess low sensitivity to perturbations, measured in terms of the
regularizer $\mathcal{V} (\Pr_n,\theta)$.

 The conclusion in Corollary \ref{cor:duality} serves as a good starting point to see why the expansion
in Theorem \ref{thm:var-reg} is plausible. Changing variables from $\lambda
\delta^{1/2}$ to $\lambda$ and $\delta^{1/2}\Delta$ to $\Delta$ in the
conclusion in Corollary \ref{cor:duality}, 
\begin{align*}
\mathcal{R}_{\delta} (\Pr_n,\theta) = \inf_{\lambda \geq 0}~\lambda
\delta^{1/2} + \frac{1}{n}\sum_{i=1}^n \sup_{\Delta} \left\{ \ell(X_i +
\delta^{1/2}\Delta,\theta) - \lambda \delta^{1/2} \Vert \Delta \Vert_q^2
\right\}.
\end{align*}
Replacing the terms $\ell(X_i + \delta^{1/2} \Delta,\theta)$ by their
respective first-order Taylor approximation, the inner suprema evaluate to 
\begin{align*}
\sup_{\Delta} \left\{\ell(X_i,\theta) + \delta^{1/2}
D_x\ell(X_i,\theta)^\intercal  \Delta - \lambda \delta^{1/2} \Vert
\Delta \Vert_q^2 \right\} = \ell(X_i,\theta) + \frac{ \delta^{1/2}}{4\lambda}%
\Vert D_x \ell(X_i,\theta)\Vert_p^2.
\end{align*}
This leads to 
\begin{align*}
\mathcal{R}_{\delta} (\Pr_n,\theta) &\approx E _{\Pr_n}
[\ell(X_i,\theta)] + \delta^{1/2} \inf_{\lambda > 0} ~\left\{\lambda + \frac{%
E _{\Pr_n}\Vert D_x \ell(X_i,\theta)\Vert_p^2}{4\lambda} \right\} \\
& = E _{\Pr_n} [\ell(X_i,\theta)] + \delta^{1/2} \mathcal{V} %
_p(\Pr_n,\theta),
\end{align*}
which heuristically justifies the conclusion in Theorem \ref{thm:var-reg}. 

Unlike the exact regularization terms exhibited in specific instances in
Examples \ref{eg:wc-return} and \ref{eg:lin-reg}, the asymptotic
regularizing effect exhibited in Theorem~\ref{thm:var-reg} holds more
broadly. More interesting is the observation that the regularization is
adapted to the model informed by the loss $\ell$, such that the resulting
sensitivity to perturbations in samples is small. Regularization terms
involving this flavour have emerged useful in adversarial training in
machine learning (see, for example, \cite{ref:goodfellow2015explaining,SHAHAM2018195, ref:volpi2018generalizing, Roth}). A principled approach towards guaranteeing adversarial robustness in machine learning contexts using Wasserstein DRO solutions has been considered in \cite{sinha2018certifiable}.


\subsection{Limiting behavior of the DRO estimator}

The most common approach towards examining the statistical properties of an
estimator is to study its limiting behavior as the number of samples grows
to infinity. While based on our discussion in the previous section we know that $\gamma=1$ should be considered, here we examine the joint limiting behavior of the triplet 
\begin{align*}
\big( \hat{\theta}_n^{\,\textnormal{erm}} , \hat{\theta}_n^{\,%
\textnormal{dro}} , \Lambda_{\delta}(\Pr_n) \big)
\end{align*}
for arbitrary $\gamma$. We must keep in mind that $\delta$ depends on $n$ but we omit this dependence to simplify the notation. 
While $%
\hat{\theta}_n^{\,\textnormal{erm}} , \hat{\theta}_n^{\,\textnormal{dro}} $ are $\mathbb{R} ^d$-valued, the third component in the
triplet, namely, $\Lambda_\delta(\Pr_n) \subseteq \mathbb{R}^d$, is the
collection of optimizers compatible with the distributional ambiguity $%
\mathcal{U}_\delta(\Pr_n)$. \textit{The limiting behavior of $\Lambda_\delta(\Pr_n)$
will help us reinforce the optimality of the projection-based prescription $%
\delta = n^{-1} \times \eta_{_{1-\alpha}}$ and the construction of DRO-based confidence regions} that are useful from the viewpoint of uncertainty quantification.
To utilize the well-known machinery of set convergence (see, for example, 
\cite{rockafellar2009variational,molchanov2005theory}) for this purpose, we
consider the right-continuous version of $\Lambda_\delta(\Pr_n)$, namely, 
\begin{align*}
\Lambda_\delta^+(\Pr_n) = \text{cl}\left(\cap_{\varepsilon > 0}
\Lambda_{\delta + \varepsilon} (\Pr_n) \right),
\end{align*}
that contains the compatible $\Lambda_\delta(\Pr_n)$ and remains closed. We
undertake this study of the limiting behavior assuming the existence of a
unique minimizer $\theta_\ast$ for \eqref{prm}, as indicated in Assumption~\ref{assume:unique-min} below.

\begin{assumption}
For each $x \in \mathbb{R} ^m$, $\ell(x,\cdot)$ is convex. Letting $%
h(x,\theta) = D_\theta \ell(x,\theta)$, there exists $\theta_\ast \in
\Theta^\circ$ satisfying the optimality conditions $E %
_{\Pr_\ast}[h(X,\theta)] = 0$, and $C \coloneqq E_{\Pr_\ast}[D_\theta
h(X,\theta_\ast)] \succ 0$. \label{assume:unique-min}
\end{assumption}

For $c > 0$, define 
\begin{align*}
b_c \coloneqq c^{1/2}C^{-1}D_\theta \mathcal{V} (\Pr_\ast,
\theta_\ast) \quad \text{ and } \quad \Lambda_c \coloneqq \left\{u:
\varphi^\ast(C u,\theta_\ast) \leq c \right\},
\end{align*}
where $D_\theta$ is the derivative operator with respect to the parameter $\theta$. Ultimately, the following result will lead us to the optimal choice $c_\ast=\eta_{_{1-\alpha}}$.

\begin{theorem}[Limit behavior]
Suppose that Assumptions \ref{assume:cost} - \ref{assume:unique-min} hold
and the distribution $\Pr_\ast$ is non-degenerate in the sense that $\Pr_\ast
(D_x\ell(X,\theta_\ast) \neq 0 ) > 0$ and that $E_{\Pr_\ast} [D_x
h(X,\theta_\ast)D_x h(X,\theta_\ast)^\intercal ] \succ 0$. Let $Z =
C^{-1}H$, where $H \sim \mathcal{N}(0, \Sigma)$ with covariance matrix $\Sigma \coloneqq 
\mathrm{Cov}_{\Pr_\ast}[h(X,\theta_\ast)]$. Then as $n \rightarrow
\infty$, we have

\begin{itemize}
\item[(i)] for $\delta = c n^{-1}$,  
\begin{align*}
\left( n^{1/2}[\hat{\theta}_n^{\,\textnormal{erm}}  - \theta_\ast], \
n^{1/2}[\hat{\theta}_n^{\,\textnormal{dro}}  - \theta_\ast], \
n^{1/2}[\Lambda_\delta^+(\Pr_n) - \hat{\theta}_n^{\,\textnormal{erm}} %
]\right) \overset{dist.}{\longrightarrow}  \big(Z, \ Z - b_c, \
\Lambda_c\big);
\end{align*}

\item[(ii)] for $\delta = c n^{-\gamma}$ with $\gamma < 1$,  
\begin{align*}
\left( n^{1/2}[\hat{\theta}_n^{\,\textnormal{erm}}  - \theta_\ast], \
n^{\gamma/2}[\hat{\theta}_n^{\,\textnormal{dro}}  - \theta_\ast], \
n^{1/2}[\Lambda_\delta^+(\Pr_n) - \hat{\theta}_n^{\,\textnormal{erm}} %
]\right) \overset{dist.}{\longrightarrow}  \big(Z, \ Z - b_c, \ 
\mathbb{R} ^d\big);
\end{align*}

\item[(iii)] and for $\delta = c n^{-\gamma}$ with $\gamma > 1$,  
\begin{align*}
\left( n^{1/2}[\hat{\theta}_n^{\,\textnormal{erm}}  - \theta_\ast], \
n^{1/2}[\hat{\theta}_n^{\,\textnormal{dro}}  - \theta_\ast], \
n^{1/2}[\Lambda_\delta^+(\Pr_n) - \hat{\theta}_n^{\,\textnormal{erm}} %
]\right) \overset{dist.}{\longrightarrow}  \big(Z, \ Z, \ \{ 0 \}\big)%
.
\end{align*}
\end{itemize}

\label{thm:joint-lim}
\end{theorem}

\noindent A proof of Theorem \ref{thm:joint-lim} is available in \cite%
{ref:blanchet2021confidence}. We focus on understanding its implications. The case where the radius of the Wasserstein ball is shrunk slower than the
recommended rate $\delta \propto n^{-1}$ (the case where $\gamma < 1$)
results in 
\begin{align*}
\Vert \hat{\theta}_n^{\,\textnormal{dro}}  - \theta_\ast \Vert =
O_p(n^{-\gamma/2}),
\end{align*}
which is suboptimal in the large-sample regime, when compared with the
benchmark set by empirical risk minimization. The rate is grossly inferior
for the choice $\delta \propto n^{-2/d}$, recommended by the use of
concentration inequalities. This characterization verifies the earlier
assertion that halving the estimation error in the DRO solution will require $2^d$ more samples. The accompanying compatible optimal solution set $%
\Lambda_\delta^+(\Pr_n) \approx \hat{\theta}_n^{\,\textnormal{erm}}  %
+ n^{-1/2}\mathbb{R} ^d$, while including the true optimum $\theta_\ast,
$ is too large to be useful in any meaningful way.

When the radius of the Wasserstein ball is shrunk at the recommended $\delta
\propto n^{-1}$, the characterization that 
\begin{align*}
\hat{\theta}_n^{\,\textnormal{dro}}  &= \theta_\ast + n^{-1/2} (Z -
b_c) + o_p(n^{-1/2}) \\
&= \hat{\theta}_n^{\,\textnormal{erm}}  - n^{-1/2} c^{1/2} C^{-1}
D_\theta \mathcal{V} (\Pr_n,\hat{\theta}_n^{\,\textnormal{erm}} %
) + o_p(n^{-1/2}),
\end{align*}
indicates the presence of an additional bias term $b_c = c^{1/2} C^{-1}
D_\theta \mathcal{V} (\Pr_\ast, \theta_\ast)$. Since $\mathcal{V} %
(\cdot)$ in \eqref{sens-defn} can be understood as the measure of
variation (or) sensitivity in the loss with respect to perturbations to
realizations of $X$, the effect of the bias term can be understood as a
``nudge" towards favouring solutions with lower sensitivity or variation
measured by $\mathcal{V} (\Pr_\ast,\theta)$. This is in line with the
adaptive regularization interpretation developed in Section \ref%
{sec:adapt-reg} for the DRO estimation \eqref{dro}.

On the other hand, for the case where the radius is shrunk faster than the
recommended rate (that is, when $\gamma > 1$), 
\begin{align*}
\hat{\theta}_n^{\,\textnormal{dro}}  = \hat{\theta}_n^{\,\textnormal{erm}}  + o_p(n^{-1/2}) \quad \text{ and } \quad
\Lambda_\delta^+(\Pr_n) = \{\hat{\theta}_n^{\,\textnormal{erm}} \} +
o_p(n^{-1/2}),
\end{align*}
revealing that there is no appreciable effect seen both in the DRO estimator 
$\hat{\theta}_n^{\,\textnormal{dro}} $ and the compatible set of
optimal solutions.

While the above discussion justifies the rate of shrinking in $\delta =
cn^{-1}$, for some $c > 0$, the optimality of the particular prescription $c
= \eta_{_{1-\alpha}}$ can be inferred from the limiting characterization as
follows. Fix $\delta = c n^{-1}$. Then $\hat{\theta}_n^{\,\textnormal{erm}} %
 - \theta_\ast = Z + o_p(n^{-1/2})$ and 
\begin{align*}
\Lambda_{\delta}^+(\Pr_n) &= \hat{\theta}_n^{\,\textnormal{erm}}  +
\Lambda_c + o_p(n^{-1/2}) \\
&= \theta_\ast + \{u + Z: \varphi^\ast(Cu,\theta_\ast) \leq c\} +
o_p(n^{-1/2}).
\end{align*}
Therefore,
\begin{align*}
\theta_\ast \in \Lambda_\delta^+ (\Pr_n) \quad \text{ if
and only if } \quad 0 \in \left\{u:\varphi^\ast(C(u-Z),\theta_\ast) \leq
c\right\} + o_p(n^{-1/2}).
\end{align*}
The optimal choice which ensures $\Lambda_\delta^+ (\Pr_n)$ is a $(1-\alpha)$%
-confidence region of $\theta_\ast$ is then given by, 
\begin{align*}
c_\ast &= \inf\{ c > 0: \theta_\ast \in \Lambda_\delta^+(\Pr_n) \text{ with
confidence } 1-\alpha\} \\
&= \inf\{ c > 0: \text{Pr}(\varphi^\ast(-CZ,\theta_\ast) \leq c) \geq
1-\alpha\} \\
&= \eta_{_{1-\alpha}},
\end{align*}
since $\eta_{_{1-\alpha}}$ is defined in Section \ref{sec:WPF-Intro-Dev} as
the $(1-\alpha)$-quantile of the limiting variable $\varphi^\ast(H,\theta_\ast)$ 
and Law($H$) = Law($-CZ$).

\subsection{Construction of DRO-compatible Confidence Regions.}
To study the confidence regions for the purposes of uncertainty quantification, we further pursue the notion of compatible confidence regions introduced in Section \ref{sec:intro-cc}. Recall that the set  $\Lambda_\delta(\Pr_n)$ serves as a projection of the distributional uncertainty $\mathcal{U}_\delta(\Pr_n)$ onto the decision space and hence can be thought of as a confidence region compatible with the DRO formulation~\eqref{dro} when suitable statistical coverage is satisfied. 

We begin this study with a Nash equilibrium characterization of the DRO estimator $\drmin.$  This can be achieved thanks to the following Theorem \ref{thm:inf-sup-exchg}.

\begin{theorem}[Inf-Sup interchange and Nash equilibrium]
Under Assumptions \ref{assume:cost} and \ref{assume:unique-min}, we have
\begin{align*}
\inf_{\theta \in \Theta} \sup_{\Pr \in \mathcal{U}_\delta(\Pr_n)} E %
_{\Pr} \left[ \ell(X,\theta)\right] = \sup_{\Pr \in \mathcal{U}%
_\delta(\Pr_n)} \inf_{\theta \in \Theta} E _{\Pr} \left[
\ell(X,\theta)\right]
\end{align*}
for any $\delta > 0$. Further, there exists a distributionally robust optimal choice $\hat{%
\theta}_n^{\,\textnormal{dro}}  \in \Lambda_\delta(\Pr_n)$ and for a
given $\varepsilon > 0$, there exists a measure $\Pr^\ast_{\varepsilon} \in 
\mathcal{U}_\delta(\Pr_n)$ such that 
\begin{align*}
\mathcal{R} (\Pr,\hat{\theta}_n^{\,\textnormal{dro}} ) -
\varepsilon \leq \mathcal{R} (\Pr^\ast_{\varepsilon},\hat{\theta}%
_n^{\,\textnormal{dro}} ) \leq \mathcal{R} (\Pr^\ast_{%
\varepsilon},\theta) + \varepsilon, 
\end{align*}
for all $\Pr \in \mathcal{U}_\delta (\Pr_n)$ and $\theta \in \Theta$. \label{thm:inf-sup-exchg}
\end{theorem}

Equipped with Theorem~\ref{thm:inf-sup-exchg}, one may view the DRO estimator $\hat{\theta}_n^{\,\textnormal{dro}} %
$ as the optimizer's strategy in a Nash equilibrium-type behavior
formed with the adversarial perturbations chosen by nature. More
importantly, the conclusion that $\hat{\theta}_n^{\,\textnormal{dro}} %
 \in \Lambda_\delta(\Pr_n)$ allows us to view the set 
\begin{align}
\Lambda_{n^{-1} \eta_{1-\alpha}}^+(\Pr_n) = \hat{\theta}_n^{\,\textnormal{erm}}  + \Lambda_{\eta_{1-\alpha}} + o_p(n^{-1/2})  \label{cr-true}
\end{align}
as a $(1-\alpha)$-confidence region simultaneously containing $\hat{\theta}%
_n^{\,\textnormal{dro}} ,\hat{\theta}_n^{\,\textnormal{erm}} $
and the unknown optimal $\theta_\ast$. A proof of Theorem \ref%
{thm:inf-sup-exchg} is available in \cite[Appendix D]%
{ref:blanchet2021confidence}.

The following characterization of the set $\Lambda_{\eta_{1-\alpha}}$, in
terms of its support function, serves as a useful tool from an algorithmic
viewpoint of constructing confidence regions. Let 
\begin{align}
\hat{\Lambda} = \bigcap_{u \in \mathbb{R} ^d}\left\{\hat{\theta}_n^{\,%
\textnormal{erm}}  + n^{-1/2}v: u^\intercal  v \leq 2[\hat{\eta%
}_{_{1-\alpha}}\hat{\varphi}(\hat{C}^{-1}u, \hat{\theta})]^{1/2} \right\}
\label{CR-supp-fn}
\end{align}
be defined in terms of any consistent estimator $\hat{\theta}$ for $\theta$, 
$\hat{C}$ for the Hessian $C$, quantile estimate $\hat{\eta}_{_{1-\alpha}}$ 
for $\eta_{_{1-\alpha}}$, and 
\begin{align}  \label{phi-savg}
\hat{\varphi}(\xi,\theta) \coloneqq 4^{-1}E _{\Pr_n} \Vert
\xi^\intercal  D_x h(X,\theta)\Vert_p^2.
\end{align}%
To see why replacing the true confidence region in the LHS of \eqref{cr-true}
by $\hat{\Lambda}$ in \eqref{CR-supp-fn} still results in an asymptotically
valid $(1-\alpha)$-confidence region, observe that the support function
of the convex set $\Lambda_c \coloneqq \{u:\varphi^\ast(Cu,\theta_\ast) \leq
c\}$ is given by
\begin{align*}
h_{\Lambda_c}(v) \coloneqq 2\left[ c\varphi(C^{-1}v)\right]^{1/2}.
\end{align*}
Then it follows from the definition of the support function that $\Lambda_c
= \cap_{u} \{v: u^\intercal  v \leq h_{\Lambda_c}(u)\}$. See~\cite{ref:blanchet2021confidence} for a more detailed explanation, examples, and
a proof of Proposition~\ref{prop:CR-alg-char} below.

\begin{proposition}[Confidence region characterization]
Under the assumptions in Theorem \ref{thm:joint-lim}, we have
\begin{align}
\textnormal{Pr}\left( \{ \theta_\ast,\hat{\theta}_n^{\,\textnormal{erm}} %
,\hat{\theta}_n^{\,\textnormal{dro}} \} \subset \hat{\Lambda}
\right) \geq 1 - \alpha + o(1) \quad \text{ as } n \rightarrow \infty.
\label{cr-guarantee}
\end{align}
\label{prop:CR-alg-char}
\end{proposition}

Thus, collecting these observations, we arrive at the recipe in Algorithm~\ref{algo:recipe} for statistical output
analysis of the DRO estimation \eqref{dro}. Note that Algorithm ~\ref{algo:recipe} outputs a confidence region which is asymptotically tight for the the parameter of interest at the prescribed $1-\alpha$ confidence level. There are many confidence regions which can be chosen, just as there are many confidence intervals that are tight at a given confidence level in the one dimensional setting. Most of the time one selects a symmetric confidence interval around the parameter of interest. However, \textit{over}estimation of a parameter of interest might be less desirable than \textit{under}estimation. As a consequence, depending on the decision maker's risk attitude, an optimal confidence interval may introduce non-symmetric features. The situation is more complex in a multi-dimensional setting. The confidence region which maximizes the likelihood (under mild assumptions) is the one that minimizes volume. Such region can be obtained by choosing a well-chosen squared Mahalanobis metric as the cost function in the Wasserstein DRO formulation, see \cite{ref:blanchet2021confidence}. However, the Wasserstein DRO formulation allows to capture different geometries which may better reflect the risk sensitivity to mis-estimation of the optimal decision. The Wasserstein DRO formulation itself informs the risk sensitivity of the modeler to such mis-estimation in connection to the loss. The confidence region obtained by Algorithm~\ref{algo:recipe} is just a natural consequence of using the Wasserstein DRO formulation to measure the impact of decision mis-estimation; see \cite{ref:blanchet2021confidence} for an illustration of the geometry induced by different Wasserstein metrics in the various optimal confidence regions obtained.

\begin{algorithm}[h!]
  \caption{A recipe for construction of a $(1-\alpha)$-confidence region satisfying \eqref{cr-guarantee}.}
  \begin{algorithmic}
  \STATE{\textbf{Input:} Training samples $X_1,\ldots,X_n$ from an unknown distribution $\Pr_\ast$.}
  \STATE{\textbf{Outline of steps:}} 
  \STATE{1. Execute Steps 1 - 2 of Algorithm~\ref{algo:recipe-dro}.} \\
  \STATE{2. Let $\hat{C}$ be the sample covariance of the collection $\{D_\theta h(X_i,\ermin):i=1,\ldots,n\}$.}\\
  \STATE{3. The set $\hat{\Lambda} = \{\ermin + n^{-1/2}v: \hat{\varphi}^\ast(\hat{C}v,\ermin) \leq \hat{\eta}_{1-\alpha}\}$ constitutes an asymptotic $(1-\alpha)$-confidence region. In the absence of the knowledge of the conjugate $\hat{\varphi}^\ast(\cdot)$, we instead compute 
\begin{align}
    \hat{\Lambda}^{(k)} = \bigcap_{i = 1, \ldots, k}\left\{\ermin + n^{-1/2}v: u_i^\T v \leq 2[\hat{\eta}_{_{1-\alpha}}\hat{\varphi}(\hat{C}^{-1}u_i, \hat{\theta})]^{1/2}  \right\},
    \label{sf-approx-cr}
\end{align}
for some large $k$, where $\{u_1, \ldots, u_k\}$ is a collection of points on $\R^d$ drawn randomly.}\\
  \STATE{4. Return the set  $\hat{\Lambda}$ or $\hat{\Lambda}^{(k)}$ as a $(1-\alpha)$-confidence region satisfying \eqref{cr-guarantee}.}
 \end{algorithmic}
    \label{algo:recipe}
    \end{algorithm}

 Due to the equivalent characterization in \eqref{CR-supp-fn}, the confidence
region $\hat{\Lambda}$ is contained in the set $\hat{\Lambda}^{(k)}$ in \eqref{sf-approx-cr} for any $k \geq 1$. Hence, albeit being larger than 
$\hat{\Lambda},$ the set $\hat{\Lambda}^{(k)}$ also serves as a $(1-\alpha)$-confidence region whose
diameter shrinks at the correct rate $O(n^{-1/2}).$ One may
replace $\hat{\theta}_n^{\,\textnormal{erm}} $ with any consistent
estimator for $\theta_\ast$ in the estimation in Steps 1 - 3  without affecting the rate of convergence.

\begin{example}[Confidence region for distributionally robust linear regression] \sloppy{Continuing the} distributionally robust linear regression
estimation in Examples \ref{eg:lin-reg} and \ref{eg:lin-reg-delta}, we have $\hat{\varphi}^\ast(z,\theta) = z^\intercal  \big( \hat{\sigma}^2 
\mathbb{I}_d + \hat{\Xi} \Vert \theta \Vert_2^2 \big)^{-1} z$ in the case $q
= 2$. Here $\hat{\sigma}$ and $\hat{\Xi}$, respectively, denote the
estimates of the error variance $\sigma^2$ and the second moment $E[XX^\T]$. The quantile $\hat{\eta}_{1-\alpha}$ is estimated as the $(1-\alpha)$-quantile of $\hat{\varphi}^\ast(H,\hat{\theta}_n^{\,\textnormal{erm}} ) = H^\intercal  \big( \hat{\sigma}^2 \mathbb{I}_d + 
\hat{\Xi} \Vert \theta \Vert_2^2 \big)^{-1} H$, where $H \sim \mathcal{N}(0,%
\hat{\Sigma})$ with $\hat{\Sigma} = \hat{\sigma}^2 \hat{\Xi}$. In this case,
we obtain $\hat{\theta}_n^{\,\textnormal{dro}} $ as a solution to %
\eqref{sqroot-lasso} with $p = 2$ and the elliptical confidence region 
\begin{align*}
\hat{\Lambda} &= \left\{ \hat{\theta}_n^{\,\textnormal{erm}}  +
n^{-1/2}v: \hat{\varphi}^\ast (\hat{C}v,\hat{\theta}_n^{\,\textnormal{erm}} %
) \leq \hat{\eta}_{_{1-\alpha}}\right\} \\
&= \hat{\theta}_n^{\,\textnormal{erm}}  + n^{-1/2}\left\{ v:
v^\intercal  \hat{A}^{-1} v \leq \hat{\eta}_{_{1-\alpha}}\right\},
\end{align*}
where $\hat{A} \coloneqq \hat{\sigma}^2 \hat{C}^{-2} + \hat{C}^{-1} \hat{%
\Sigma}^{-1} \hat{C}^{-1} \Vert \hat{\theta}_n^{\,\textnormal{erm}}  %
\Vert_2^2$.
\label{eg:CR-lin-reg-eg} 
\end{example}

\subsection{Finite-sample error bounds}
\label{sec:fin-sample-bnds}

The goal of this subsection is to briefly discuss the elements involved in finite sample error bounds obtained in the literature. The literature considers the cases where the ambiguity radius $\delta$ is taken to be (i) non-vanishing with the sample size $n$, and (ii) when $\delta$ is decreased with nominal dependence on the ambient dimensions $m$.

The value in non-asymptotic bounds lies in the information provided in terms of various distributional parameters. This information could be helpful if there is a design parameter that can be used to mitigate a small sample size. For use of the formulation \eqref{dro} where the ambiguity radius $\delta$ is non-vanishing with the sample size, we have the following result from~\cite[Theorem 2]{ref:lee2018minimax}. For stating the result, let $\mathcal{C}$ denote the Dudley entropy integral~\cite{talagrand2014upper} for the function class $\{\ell(\cdot,\theta): \theta \in \Theta\}$. 

\begin{theorem}[Finite sample guarantee for non-vanishing $\delta$]
Suppose that $\xsp$ is a bounded subset of $\R^m$ and the collection of functions $\{\ell(\cdot,\theta): \theta \in \Theta\}$ are uniformly bounded and $L$-Lipschitz, that is, there exist positive constants $M$, $L$ such that $0 \leq \ell(x,\theta) \leq M$ and $\vert \ell(x,\theta) - \ell(x^\prime,\theta) \vert \leq L \Vert x - x^\prime\Vert$ for all  $x$, $x^\prime \in \mathcal{X}$ and $\theta \in \Theta$. Then for the transportation cost $c(x,x^\prime) = \Vert x  - x^\prime \Vert^r$, $r \in [1,\infty)$, we have
\begin{align*}
\drorisk(\Pr_\ast,\drmin) - \inf_{\theta \in \Theta} \drorisk(\Pr_\ast,\theta) \leq  n^{-1/2} \left[ c_0 + c_1 \delta^{-1+1/r} + c_3 \log(2/\varepsilon)\right],         
\end{align*}
with probability at least $1 - \varepsilon$. With $\mathcal{C}$ denoting the Dudley entropy integral for the function class $\{\ell(\cdot,\theta): \theta \in \Theta\}$, the constants $c_0$, $c_1$ and $c_2$ are identified as follows:
\begin{align*}
    c_0 \Let 48 \mathcal{C}, \quad c_1 
    \Let 48L\cdot \textnormal{diam}(\xsp)^r, \quad \text{ and } \quad c_3 \Let 2^{-1/2} \times 3M.
\end{align*}
\end{theorem}

The above finite sample bound is intended for use in settings such as domain adaptation (see \cite{ref:lee2018minimax} and references) where it is meaningful to consider non-vanishing radius $\delta$. However for instances where $\delta$ is taken to decrease with the sample size, the above finite sample bound is of relatively limited utility; 
for example, for the choice $\delta = O(1/n)$ exhibited in Section \ref{sec:dro-clt}, we have the resulting finite sample error to be $O(1)$. Refined finite sample guarantees which are applicable for vanishing choices of ambiguity radius $\delta$ have been developed in \cite{ref:shafieezadeh2017regularization,chen2018robust,gao2020finite}. Considering distributionally robust formulations of the form 
\begin{align}
     \inf_{\theta \in \Theta} \sup_{\Pr \in \mathcal{U}_\delta(\Pr_n)} \E_{\Pr} \left[ \ell(\theta^\T X, Y)\right] 
     \label{supervised-dro}
\end{align}
motivated by supervised learning problems, \cite{ref:shafieezadeh2017regularization} develops the following generalization bound. With $\{e_i: i =1,\ldots, m+1\}$ denoting the standard basis vectors in $\R^{m+1},$ let  $M_{m,p} \Let \max_{i \leq m+1} \Vert e_i \Vert_p$ for use in Theorem \ref{thm:fs-2} below. 

\begin{theorem}[Finite sample guarantee with $\delta$ decreasing in sample size]
\label{thm:fs-2}
 Let the ground transportation cost be of the form $c((x,y), (x^\prime,y^\prime)) = \Vert (x,y)  -  (x^\prime, y^\prime) \Vert_q$. Suppose $X$ is light-tailed in the sense there exist constants $a > 1$, $A > 0$ such that $\E_{\Pr_\ast}[\exp(a\Vert (X,Y) \Vert_p)] \leq A$. For $p$ satisfying $1/p + 1/q = 1$, suppose the loss $\ell(\cdot)$ and the set $\Theta$ are such that either one of the following holds:
 \begin{itemize}
     \item[(i)] $\ell(x,y) = L(\theta^\T x - y)$ for a Lipschitz $L:\R \rightarrow \R$ and $\inf_{\theta \in \Theta} \Vert (\theta, -1)\Vert_p \geq  \underline{\Omega}$  if \eqref{supervised-dro} is a regression problem; or
     \item[(ii)]  $\ell(x,y) = L(y\theta^\T x)$ for a Lipschitz $L:\R \rightarrow \R$ and $\inf_{\theta \in \Theta} \Vert \theta\Vert_p \geq  \underline{\Omega}$ if \eqref{supervised-dro} is a classification problem. 
 \end{itemize}
 Then there exist constants $c_0 \geq 1$, $c_1 > 0$ depending only on $a$, $A$ such that for any \sloppy{$n \geq \max\{ (16m/c_1)^2,16\log(c_0/\varepsilon)/c_1)\}$} and the ambiguity radius choice
 \begin{align*}
     \delta \geq \frac{2\textnormal{diam}(\Theta)}{\sqrt{n}\underline{\Omega}} \left[ m  A M_{m,p} + \left[c_1^{-1}(m\log\sqrt{n} + \log (c_0\varepsilon^{-1}))\right]^{1/2} \right],
 \end{align*}
 we have the following generalization bound holding with probability at least $1-\varepsilon:$
 \begin{align*}
     \sup_{\theta \in \Theta} ~\left[\risk(\Pr_\ast,\theta) - \drorisk(\Pr,\theta)\right] \leq 0. 
     \end{align*}
\end{theorem}
The choice of $\delta$ in the previous result almost matches (up to a logarithmic factor) the optimal choice based on the Wasserstein Profile Function, however, the guarantees are non-asymptotic (albeit, under stronger conditions). The finite sample error bounds in \cite{gao2020finite} do match the optimal decay rates and are available also for the squared cost choice $c(x,x^\prime) = \Vert x - x^\prime \Vert^2$. However, the constants involved in the bound are not immediately computable in terms of the elements in the exposition here. The assumptions involve transportation inequalities (implying super-exponentially decaying tails) to be satisfied by the underlying distributions. We refer the readers to \cite{gao2020finite} for these refined finite sample guarantees. Additional finite sample guarantees with optimal convergence rates for regression problems (under bounded support) are derived in \cite{chen2018robust}.

\section{Conclusions and Final Considerations} \label{sec:conclusions}
Data-driven Wasserstein-DRO formulations have gained significant attention in recent years due to their intuitive appeal as a direct mechanism to improve out-of-sample performance and generalization. A key ingredient in these applications is the choice of the uncertainty size (the radius $\delta$ of the ambiguity set in our discussion). Most of the Wasserstein-DRO literature either advocates a choice of $\delta$ which either suffers from the curse of dimensionality (e.g., by enforcing the underlying data-generating distribution to be in the ambiguity set) or a choice with strong (non-asymptotic) guarantees in terms of generalization bounds at the expense of difficult-to-compute constants or strong assumptions in the underlying distributions. The main practical method for choosing $\delta$ is cross validation (CV), which can be safe if used properly (following the prescriptions of \cite{ref:shao1993linear} at least in the linear regression setting). However, CV could be time consuming and very data intensive.

We have focused here on a method to choose $\delta$ that is based on the projection analysis in the Wasserstein geometry. This method provides easy-to-implement algorithmic procedures with solid statistical guarantees (including coverage and asymptotic normality). As a by-product, this method introduces a new hypothesis test with its test statistic being computed from the projection distance. The application of these types of methods in the study of optimal transport formulation for data-driven \textit{decisions} is an area of research which remains to be explored. Of significant research interest, but outside of the scope of this tutorial, is the use of related optimal transport-related DRO formulations which would avoid the curse of dimensionality in the selection of $\delta$; these include, for example, the so-called sliced-Wasserstein distance \cite{ref:kolouri2019generalized}, the smoothed Wasserstein distance \cite{ref:goldfeld2020convergence, ref:goldfeld2020asymptotic} among others. 

Under regularity conditions, the optimality gap obtained by the empirical risk minimization solution, technically defined as the difference $\mathcal{R} \big(\Pr_\ast, {\hat{\theta}_n^{\,\textnormal{erm}} } \big) - \inf_{\theta \in \Theta} \mathcal{R} (\Pr_\ast, \theta)$, is asymptotically optimal as the sample size increases in the second-order convex sense compared to a wide range of regularizations including the DRO-type formulations~\cite{ref:lam2021on}. While this optimality property is remarkable, it is important to keep in perspective the assumptions imposed therein. For example, conditions such as the existence of a unique optimizer, twice differentiability at the optimum, and a fixed dimensionality environment appear key in the development of the asymptotic optimality result in \cite{ref:lam2021on}. 
Furthermore, in more complex decision making tasks, the notion of optimality gap may require refinements. For a concrete example, a portfolio allocation is usually based on historical data, and it is natural to evaluate the optimality gap of the portfolio selection \textit{conditional} on the side information, for instance, based on market implied volatilities which reflect future market expectations. Several DRO formulations have been proposed recently to accommodate additional information; see, for example, \cite{ref:blanchet2021optimal,ref:esteban2020distributionally, kannan2020residualsbased, ref:nguyen2021robustifying}, and this remains an emerging future direction for research.

\textbf{Acknowledgements.} Material in this paper is based upon work supported by the Air Force Office of Scientific Research under award number FA9550-20-1-0397. Support from Singapore MOE-SUTD research grant SRG-ESD-2018-134 and MOE Academic Research Fund MOE2019-T2-2-163 are gratefully acknowledged. We also appreciate additional support from NSF grants 1915967, 1820942, 1838576.

\bibliographystyle{TutORials}
\bibliography{main.bbl}

\end{document}